\pdfoutput=1

\documentclass[sigplan,10pt,nonacm]{acmart}

\setcopyright{acmcopyright}
\copyrightyear{2018}
\acmYear{2018}
\acmDOI{XXXXXXX.XXXXXXX}

\acmConference[Conference acronym 'XX]{Make sure to enter the correct
  conference title from your rights confirmation emai}{June 03--05,
  2018}{Woodstock, NY}
\acmPrice{15.00}
\acmISBN{978-1-4503-XXXX-X/18/06}

\usepackage{xifthen}
\usepackage{xspace}
\usepackage{tikz-cd}

\usepackage{xurl}
\usepackage{setouterhbox}
\usepackage{etoolbox}

\makeatletter
\patchcmd\deferred@thm@head{\sbox\@labels{\normalfont#1}}{\setouterhbox\@labels \normalfont#1\endsetouterhbox}{}{\fail}
\makeatother

\usetikzlibrary{decorations.markings}
\tikzset{mid vert/.style={/utils/exec=\tikzset{every node/.append style={outer sep=0ex}},
		postaction=decorate,decoration={markings,
			mark=at position 0.5 with {\draw[-] (0,#1) -- (0,-#1);}}},
	mid vert/.default=0.75ex}

\newcommand{\longhash}{9e43b0d347aea512a2f999f5e3581e8c89ae7236}
\newcommand{\shorthash}{9e43b0d}

\newcommand{\coqdocbasebaseurl}{https://benediktahrens.gitlab.io/unimathdoc3/\shorthash}

\newcommand{\coqdocbaseurl}{\coqdocbasebaseurl /UniMath.}
\newcommand{\urlhash}{\#}
\newcommand{\coqdocurl}[2]{\coqdocbaseurl #1.html\urlhash #2}
\newcommand{\nolinkcoqident}[1]{\color{blue}{\nolinkurl{#1}}} %
\makeatletter
\newcommand{\coqident}{\begingroup\@makeother\#\@coqident}
\newcommand{\@coqident}[3][]{%
  \ifthenelse{\isempty{#2}}%
  {\nolinkcoqident{#3}}%
  {\ifthenelse{\isempty{#1}}%
  {\href{\coqdocurl{#2}{#3}}{\nolinkcoqident{#3}}}%
  {\href{\coqdocurl{#2}{#3}}{\nolinkcoqident{#1}}}}%
\endgroup}
\newcommand{\coqfile}[2]{%
  \ifthenelse{\isempty{#1}}%
  {\href{\coqdocbaseurl #2.html}{\nolinkcoqident{#2.v}}}%
  {\href{\coqdocbaseurl #1.#2.html}{\nolinkcoqident{#2.v}}}}
\makeatother
\newcommand{\commit}{\href{https://github.com/UniMath/UniMath/tree/\longhash}{\shorthash}}

\AtEndPreamble{%
\theoremstyle{acmdefinition}
\newtheorem{constrInternal}[theorem]{Construction}%
\newtheorem{remark}[theorem]{Remark}
}

\usepackage[draft]{optional}
\newcommand{\plan}[1]{}
\newcommand{\BA}[1]{}
\newcommand{\PN}[1]{}
\newcommand{\NR}[1]{}
\newcommand{\NW}[1]{}
\opt{draft}{
\renewcommand{\plan}[1]{\textcolor{blue}{TODO: #1}\PackageWarning{TODO}{TODO: #1}}
\renewcommand{\BA}[1]{\textcolor{orange}{BA: #1}\PackageWarning{TODO}{TODO: #1}}
\renewcommand{\PN}[1]{\textcolor{purple}{PN: #1}\PackageWarning{TODO}{TODO: #1}}
\renewcommand{\NR}[1]{\textcolor{teal}{NR: #1}\PackageWarning{TODO}{TODO: #1}}
\renewcommand{\NW}[1]{\textcolor{magenta}{NW: #1}\PackageWarning{TODO}{TODO: #1}}
}

\newcommand{\teletype}[1]{\ensuremath{\mathtt{#1}}}
\newcommand{\systemname}[1]{\teletype{\color{darkgray}#1}\xspace}
\newcommand{\UniMath}{\systemname{UniMath}}
\newcommand{\Coq}{\systemname{Coq}}

\newcounter{saveenumi}
\newcommand{\saveitem}{\setcounter{saveenumi}{\value{enumi}}}
\newcommand{\restoreitem}{\setcounter{enumi}{\value{saveenumi}}}

\usepackage[capitalize]{cleveref}

\newcommand{\constfont}[1]{\ensuremath{\mathsf{#1}}}
\newcommand{\cat}[1]{\ensuremath{\constfont{#1}}\xspace}
\newcommand{\disp}[1]{\overline{#1}} %

\newcommand{\depEq}[0]{=_*} %

\newcommand{\CC}[0]{\cat{C}} %
\newcommand{\CD}[0]{\cat{D}} %
\newcommand{\SET}[0]{\cat{Set}} %
\newcommand{\id}[1]{\cat{id}_{#1}} %
\newcommand{\idI}[0]{\cat{id}} %
\newcommand{\inv}[1]{#1^{-1}} %
\newcommand{\nattrans}[2]{#1 \Rightarrow #2} %
\newcommand{\iso}[2]{#1 \cong #2} %

\newcommand{\CB}[0]{\cat{B}} %
\newcommand{\DBob}[2]{#1_{#2}} %
\newcommand{\DBmor}[3]{#1 \rightarrowD{#3} #2} %
\newcommand{\DBcell}[3]{#1 \Rightarrow_{#3} #2} %
\newcommand{\cell}[2]{#1 \Rightarrow #2} %
\newcommand{\adjequiv}[2]{#1 \simeq #2} %
\newcommand{\invcell}[2]{#1 \cong #2} %

\newcommand{\midhor}{
	\hspace{-0.5em}
  \begin{tikzcd}[sep=0.2in]
		 \arrow[r, mid vert] \pgfmatrixnextcell \
	\end{tikzcd}
  \hspace{-.7em}
}
\newcommand{\hor}[2]{{#1} \midhor {#2}} %
\newcommand{\horid}[1]{\cat{id}_{#1}} %
\newcommand{\horcomp}[2]{#1 \odot #2} %
\renewcommand{\rightarrow}{
  \hspace{-0.5em}
  \begin{tikzcd}[sep=0.2in]
		 \arrow[r] \pgfmatrixnextcell \
	\end{tikzcd}
  \hspace{-.7em}
}
\renewcommand{\to}{
  \hspace{-0.5em}
  \begin{tikzcd}[sep=0.2in]
		 \arrow[r] \pgfmatrixnextcell \
	\end{tikzcd}
  \hspace{-.7em}
}
\newcommand{\rightarrowD}[1]{
  \hspace{-0.5em}
  \begin{tikzcd}[sep=0.2in]
		 \arrow[r] \pgfmatrixnextcell \! \! \! _{#1} \ \
	\end{tikzcd}
  \hspace{-.7em}
}
\newcommand{\Dsquare}[4]{\left( #1 \, \, \, ^{#3}_{#4} \, \, \, #2 \right)} %
\newcommand{\sqvid}[1]{\cat{id}^{v}_{\cat{sq}}(#1)} %
\newcommand{\sqvcomp}[2]{#1 \cdot_{\cat{sq}} #2} %
\newcommand{\sqhid}[1]{\cat{id}^{h}_{\cat{sq}}(#1)} %
\newcommand{\sqhcomp}[2]{#1 \odot_{\cat{sq}} #2} %
\newcommand{\Dlunitor}[1]{\lambda_{#1}}
\newcommand{\Drunitor}[1]{\rho_{#1}}
\newcommand{\Dassociator}[3]{\alpha_{(#1, #2, #3)}}

\newcommand{\Dob}[3]{#1_{(#2, #3)}} %
\newcommand{\Dmor}[4]{#1 \rightarrowD{(#3, #4)} #2} %
\newcommand{\Did}[1]{\disp{\cat{id}}_{#1}} %
\newcommand{\total}[1]{\int #1} %
\newcommand{\projl}[1]{\pi^{#1}_1} %
\newcommand{\projr}[1]{\pi^{#1}_2} %
\newcommand{\idtoiso}[2]{\cat{idtoiso}_{#1, #2}} %

\newcommand{\Arr}[1]{\cat{Arr}(#1)}
\newcommand{\Comma}[2]{\cat{Comma}(#1, #2)}
\newcommand{\Span}[1]{\cat{Span}(#1)}
\newcommand{\StructCospan}[1]{\cat{StructCospan}(#1)}
\newcommand{\Lens}[1]{\cat{Lens}(#1)}

\newcommand{\SpanDiag}[5]{#1 \xleftarrow{#2} #3 \xrightarrow{#4} #5} %
\newcommand{\Cospan}[5]{#1 \xrightarrow{#2} #3 \xleftarrow{#4} #5} %

\newcommand{\proddispbicat}[2]{#1 \times #2} %
\newcommand{\fullsubbicat}[2]{} %

\newcommand{\bicatofcats}[0]{\cat{UnivCat}} %

\newcommand{\dispbicatoftwosideddispcat}[0]{\cat{TwoSidedDisp}_{\cat{d}}} %
\newcommand{\bicatoftwosideddispcat}[0]{\cat{TwoSidedDisp}} %

\newcommand{\dispbicathorid}[0]{\cat{HorId}_{\cat{d}}} %
\newcommand{\dispbicathorcomp}[0]{\cat{HorComp}_{\cat{d}}} %
\newcommand{\dispbicathoridcomp}[0]{\cat{HorIdComp}_{\cat{d}}} %
\newcommand{\bicathoridcomp}[0]{\cat{HorIdComp}} %

\newcommand{\dispbicatlunitor}[0]{\cat{Lun}_{\cat{d}}} %
\newcommand{\dispbicatrunitor}[0]{\cat{Run}_{\cat{d}}} %
\newcommand{\dispbicatassociator}[0]{\cat{Assoc}_{\cat{d}}} %
\newcommand{\dispbicatunitorassociator}[0]{\cat{UnAssoc}_{\cat{d}}} %
\newcommand{\bicatunitorassociator}[0]{\cat{UnAssoc}} %

\newcommand{\bicatofdoublecats}[0]{\cat{DoubleCat}} %

\newcommand{\doublecatvercat}[1]{} %
\newcommand{\doublecathor}[1]{} %
\newcommand{\doublecathorid}[2]{\cat{id}_{\cat{h}}^{#1}(#2)} %
\newcommand{\doublecathoridsq}[2]{\cat{id}_{\cat{sq}}^{#1}(#2)} %
\newcommand{\doublecathorcomp}[2]{#1 \odot #2} %
\newcommand{\doublecatlunitor}[1]{} %
\newcommand{\doublecatrunitor}[1]{} %
\newcommand{\doublecatassociator}[1]{} %

\newcommand{\doublefunctorhorid}[2]{#1_{\cat{id}}(#2)} %
\newcommand{\doublefunctorhorcomp}[3]{#1_{\cat{comp}}(#2, #3)} %

\newcommand{\lensput}[1]{\cat{put}_{#1}}
\newcommand{\lensget}[1]{\cat{get}_{#1}}

\begin{document}

\title{Univalent Double Categories}

\author{Niels van der Weide}
\email{nweide@cs.ru.nl}
\orcid{0000-0003-1146-4161}
\affiliation{%
  \institution{Radboud University}
  \city{Nijmegen}
  \country{The Netherlands}
}
\author{Nima Rasekh}
\email{rasekh@mpim-bonn.mpg.de}
\orcid{0000-0003-0766-2755}
\affiliation{%
  \institution{Max Planck Institute for Mathematics}
  \city{Bonn}
  \country{Germany}
}

\author{Benedikt Ahrens}
\email{B.P.Ahrens@tudelft.nl}
\orcid{0000-0002-6786-4538}
\affiliation{%
  \institution{Delft University of Technology}
  \country{The Netherlands}
}
\affiliation{%
  \institution{University of Birmingham}
  \country{United Kingdom}
}

\author{Paige Randall North}
\email{p.r.north@uu.nl}
\orcid{0000-0001-7876-0956}
\affiliation{%
  \institution{Utrecht University}
  \country{The Netherlands}
}

\renewcommand{\shortauthors}{Van der Weide, Rasekh, Ahrens, and North}

\begin{abstract}
  Category theory is a branch of mathematics that provides a formal framework for understanding the relationship between mathematical structures.
  To this end, a category not only incorporates the data of the desired objects, but also ``morphisms", which capture how different objects interact with each other.
  Category theory has found many applications in mathematics and in computer science, for example in functional programming.
  
  Double categories are a natural generalization of categories which incorporate the data of two separate classes of morphisms, allowing a more nuanced representation of relationships and interactions between objects. Similar to category theory, double categories have been successfully applied to various situations in mathematics and computer science, in which objects naturally exhibit two types of morphisms. Examples include categories themselves, but also lenses, petri nets, and spans.
  
  While categories have already been formalized in a variety of proof assistants, double categories have received far less attention. In this paper we remedy this situation by presenting a formalization of double categories via the proof assistant Coq, relying on the Coq UniMath library. As part of this work we present two equivalent formalizations of the definition of a double category, an unfolded explicit definition and a second definition which exhibits excellent formal properties via 2-sided displayed categories. As an application of the formal approach we establish a notion of univalent double category along with a univalence principle: equivalences of univalent double categories coincide with their identities. 
\end{abstract}



\received{20 February 2007}
\received[revised]{12 March 2009}
\received[accepted]{5 June 2009}

\maketitle

\section{Introduction}

Double categories \cite{ehresmann1963categories} are a categorical concept that captures more structure than a category. They are often succinctly defined as an internal (pseudo)category in the 2-category of categories.
A double category has objects, two kinds of morphisms --- called \emph{vertical} and \emph{horizontal}, respectively --- and fillers for squares formed from horizontal and vertical morphisms.
As such, a double category can capture two different kinds of morphisms (and their interplay) between mathematical objects.

Many mathematical objects are better understood within a double category than within a category;
for instance, the double category of sets, functions, and relations. The objects of this double category are sets $X$, the vertical morphisms $X \to Y$ are functions $X \to Y$, the horizontal morphisms $\hor{X} Y$ are relations, i.e. subsets of $X \times Y$. Considering this double category allows one to generalize classical set theory (largely overlapping with the generalization given by topos theory). Similarly, one can also consider the double category of categories, functors, and profunctors, and this has been used to great success to generalize category theory \cite{MR2844536}.

Applications of double categories have become ubiquitous in mathematics and computer science; see, for instance, its applications in systems theory \cite{courser2020open,Jaz_Myers_2021,Baez2022structuredversus} and
programming languages theory \cite{DBLP:conf/lics/DagandM13,DBLP:conf/fossacs/NewL23}.

In the present work we develop the notion of univalent double category and a library of univalent double categories in univalent foundations.
Our main result states that the bicategory of univalent double categories is univalent.
As a consequence, the type of identities $A = B$ between univalent double categories $A$ and $B$ coincides with the type $A \simeq B$ of equivalences from $A$ to $B$.
The proof of this result relies crucially on Voevodsky's univalence axiom.
The result entails that any construction on univalent double categories can be transported across equivalences --- an instance of the \emph{univalence principle} \cite{up}.

Double categories consist of a lot of data (see \cref{sec:double-categ-their}),
and morphisms of double categories --- and morphisms between these morphisms --- need to preserve that structure suitably.
In other words, the bicategory of double categories is quite complicated.
For this reason, a naïve, brute-force approach to proving univalence of this bicategory would lead to difficult proofs.
Instead, we develop technology to build the bicategory of double categories in layers, using displayed bicategories \cite{DBLP:journals/mscs/AhrensFMVW21}.
We then prove every layer univalent, and obtain that their ``total bicategory'' --- which is the desired bicategory of univalent double categories --- is univalent, by a result from \cite{DBLP:journals/mscs/AhrensFMVW21}.
The key layer we consider in this approach is the layer of ``2-sided displayed categories''.
These gadgets are a simple variation on the notion of displayed category \cite{DBLP:journals/lmcs/AhrensL19}.
Through their use, we derive a modular proof of univalence of the bicategory of univalent double categories;
in particular, we can reuse an existing proof of univalence of univalent categories from \cite{DBLP:journals/mscs/AhrensFMVW21}.

Building the bicategory of univalent double categories in a layered way also gives rise to an interesting characterization of equivalences of double categories.
In \cref{sec:equiv-inv-double}, we show that a double functor between univalent double categories is an adjoint equivalence if it is a strong double functor and an adjoint equivalence on the underlying 2-sided displayed category.

\subsection{UniMath}
\label{sec:overv-univ-found}

In this section we provide a brief introduction to univalent foundations and \UniMath, and fix notations used throughout the paper.
By univalent foundations, we mean Martin-Löf type theory (MLTT) plus Voevodsky's univalence axiom.
We use standard notation for the type and term formers of MLTT;
in particular, we write $a = b$ for the type of identifications/equalities/paths from $a$ to $b$. %

Crucially, we rely on the notion of \emph{homotopy level}, and, in particular, the notions of proposition  and set  of univalent foundations:
a type $X$ is a proposition if $\prod_{x,y:X}x=y$ is inhabited,
and a set if the type $x = y$ is a proposition for all $x, y : X$.
Hence, despite working in \Coq, we do not rely on the universes \coqident{}{Prop} or \coqident{}{Set}.

We do not rely on any inductive types other than the ones specified in the prelude of \UniMath, such as identity types, sum types, natural numbers, and booleans.

Our key result, the univalence principle for univalent double categories, relies on the univalence principle for types, also known as Voevodsky's Univalence Axiom.
This axiom is added to \Coq as a postulate in \UniMath.

\subsection{Related Work}
\label{sec:related-work}

Double categories have been formalized in several computer proof assistants.

Murray, Pronk, and Szyld~\cite{doublecats_lean} worked towards defining double categories in the Lean proof assistant.
The chosen approach is to define double categories as category objects in the category of categories.
The corresponding pull request%
\footnote{\url{https://github.com/leanprover-community/mathlib/pull/18204}}
to the mathlib library seems to have been abandoned.

Hu and Carette~\cite{DBLP:conf/cpp/HuC21} started a library of category theory in Agda.
At the time of writing that article, ``[\ldots] double categories [\ldots]  are still awaiting'' formalization.
In the meantime, the definition of double categories, as well as the construction of the dual of a double category (swapping horizontal and vertical morphisms),
have been implemented.%
\footnote{\url{https://github.com/agda/agda-categories/blob/36abe6bff98be027bd4fcc3306d6dac8b2140079/src/Categories/Double/Core.agda}}

In \cref{sec:double-categ-their}, we give a more detailed comparison between the different notions of double category mentioned here.

\subsection{Computer Formalization}
\label{sec:comp-form}

The formalization accompanying this paper is based on the UniMath library \cite{UniMath}, a library of computer-checked mathematics in the univalent style.
UniMath itself is based on the Coq proof assistant.

Our code has been integrated in the UniMath library in commit \commit.
From this commit, we compiled an HTML documentation of UniMath;
throughout this article we include links to this documentation, as in \coqident{Bicategories.DisplayedBicats.DispBicat}{disp_bicat}.
The interested reader can type-check our definitions by following the compilation instructions of the UniMath library.

\subsection{Synopsis}

In \cref{sec:double-categ-their} we informally review and motivate the notion of double category, and give an elementary, unfolded definition.
The unfolded definition is easy to understand, but proving a univalence principle for it \emph{directly} would be tough.
For this reason, we introduce, in \cref{sec:2-sided-displayed}, the notion of \emph{2-sided displayed categories}; we use these in \cref{sec:bicat-of-double-cat} to build a bicategory of displayed categories that does not make use of the elementary definition.
To prepare for this construction, we review the notion of displayed bicategory in \cref{sec:disp-bicat}.
In \cref{sec:exa-double-cats} we construct several examples of univalent double categories.
In \cref{sec:equiv-inv-double} we give a characterization of adjoint equivalences, and of invertible 2-cells, in the bicategory of univalent double categories --- that is, of equivalences of double categories and of invertible transformations between functors of univalent double categories.

\section{Double Categories}
\label{sec:double-categ-their}
In this section, we give a brief overview of the theory of double categories. 
Intuitively, a double category is a category with an extra class of morphisms.
Morphisms in one class of morphisms is called \emph{vertical} morphisms,
and morphisms in the other class are called \emph{horizontal} morphisms.
We see horizontal morphisms as ``extra'' morphisms, and for those, the laws do not hold up to equality (see \cref{rem:double-cat-flavors}).
We denote the vertical morphisms by $x_1 \rightarrow x_2$ and the horizontal morphisms by $\hor{x}{y}$.
In addition, any double category features a collection of squares, parametrized by a boundary consisting of two horizontal and two vertical morphisms with ``compatible'' endpoints as follows:
\[\begin{tikzcd}
	{x_1} & {y_1} \\
	{x_2} & {y_2}
	\arrow["{v_1}"', from=1-1, to=2-1]
	\arrow["{h_1}", "\shortmid"{marking}, from=1-1, to=1-2]
	\arrow["{v_2}", from=1-2, to=2-2]
	\arrow["{h_2}"', "\shortmid"{marking}, from=2-1, to=2-2]
\end{tikzcd}\]
Such squares are also denoted as $\Dsquare{v_1}{v_2}{h_1}{h_2}$.
For both the horizontal and the vertical morphisms we have identities and compositions.
However, there is an essential difference between the two classes of morphisms:
laws for the vertical morphisms hold up to equality,
whereas the laws of horizontal morphisms hold up to a square.
Concretely, this means that we have \emph{unitor} and \emph{associator} squares that witness the unitality and associativity of horizontal composition.
In addition, this data is coherent: we also require the triangle and pentagon equation for this data.

\begin{remark}
\label{rem:double-cat-flavors}
The notion of double category comes in several flavors.
For example, there is the notion of \emph{strict double category},
and in those, unitality and associativity of composition holds as an equality.
However, in the remainder of this paper, we look at \emph{pseudo double categories}, and in those, composition of horizontal morphisms is only weakly unital and associative.
Pseudo double categories are a useful generalization of strict double categories.
Some examples, such as spans (\cref{exa:span-double-cat}) and structured cospans (\cref{exa:struct-cospan-double-cat}),
are pseudo double categories, but not strict ones.
\end{remark}

Double categories play a prominent role in applied category theory.
For example, Clarke defined a double category of lenses \cite{Clarke2023}, and lenses have become an important tool in the study of databases and datatypes; see, e.g., \cite{DBLP:conf/pods/BohannonPV06}.
In addition, Baez and Master \cite{DBLP:journals/mscs/BaezM20} defined a double category of Petri Nets, which are used in the study of parallel programs \cite{DBLP:conf/mfcs/Kotov78} and modeling hardware \cite{DBLP:journals/csur/Peterson77}.
Baez and Courser defined a double category of structured cospans and of decorated cospans \cite{MR3384097,MR4170469}, which are used to model open systems.

There are several approaches to defining the notion of double category, and each comes with their own merits and drawbacks.
The most concise definition is that a double category is a pseudocategory internal to the bicategory of categories.
While this definition is clean and short, its drawback is that composition is described using pullbacks,
which makes it more cumbersome to work with.
More concretely, let us assume we have two categories $\CC_H$ and $\CC_V$ together with functors $S, T : \CC_H \rightarrow \CC_V$.
If we were to use this definition,
then a horizontal arrow from $x : \CC_V$ to $y : \CC_V$ would consist of an object $h : \CC_H$ together with isomorphisms $\iso{S(h)}{x}$
and $\iso{T(h)}{y}$.
In addition, the composition operation for horizontal arrows would take three objects $x, y, z : \CC_V$,
two horizontal arrows $h, k : \CC_H$
and isomorphisms $\iso{S(h)}{x}$, $\iso{T(h)}{y}$, $\iso{S(k)}{y}$, and $\iso{T(k)}{z}$,
and it returns a horizontal arrow $h \cdot k : \CC_H$ together with isomorphisms $\iso{S(h \cdot k)}{x}$ and $\iso{T(h \cdot k)}{z}$.

\begin{remark}
\label{rem:strict}
Note that one could also look at categories internal to a 1-category instead of a bicategory.
By looking at categories internal to the 1-category of strict categories,
one obtains yet another notion of double category.
This approach is taken in Lean \cite{doublecats_lean,X20}, where pullbacks are used directly,
and in 1lab \cite{1lab}, where pullbacks are avoided by looking at the internal language of a presheaf category.
However, this approach comes with a significant limitation:
by looking at strict categories, one loses examples such as spans in $\SET$ (\cref{exa:span-double-cat}),
and the square construction for univalent categories (\cref{exa:square-double-cat}).
Note that if one assumes uniqueness of identity proofs,
then internal categories in the 1-category of strict categories corresponds to strict double categories as discussed in \cref{rem:double-cat-flavors}.
\end{remark}

We can avoid pullbacks by going for an unfolded definition, which looks as follows:

\begin{definition}[\coqident{Bicategories.DoubleCategories.DoubleCatsUnfolded}{doublecategory}]
\label{def:double-cat-unfolded}
A \textbf{double category} consists of
\begin{enumerate}
  \item\label{double-cat:vertical-cat} a category $\CC$ called the \textbf{vertical category};
  \item\label{double-cat:hor-mor} for all objects $x : \CC$ and $y : \CC$, a type $\hor{x}{y}$ of \textbf{horizontal morphisms};
  \item\label{double-cat:hor-id} for every object $x : \CC$ a \textbf{horizontal identity} $\horid{x} : \hor{x}{x}$;
  \item\label{double-cat:hor-comp} for all horizontal morphisms $h : \hor{x}{y}$ and $k : \hor{y}{z}$, a \textbf{horizontal composition} $\horcomp{h}{k} : \hor{x}{z}$;
  \item\label{double-cat:squares} for all horizontal morphisms $h : \hor{x_1}{y_1}$ and $k : \hor{x_2}{y_2}$ and vertical morphisms $v_x : \hor{x_1}{x_2}$ and $v_y : \hor{y_1}{y_2}$, a set $\Dsquare{v_x}{v_y}{h}{k}$ of \textbf{squares};
  \item\label{double-cat:vertical-id-square} for all horizontal morphisms $h : \hor{x}{y}$, we have a \textbf{vertical identity} $\sqvid{h} : \Dsquare{\id{x}}{\id{y}}{h}{h}$;
  \item\label{double-cat:vertical-comp-square} for all squares $\tau_1 : \Dsquare{v_1}{w_1}{h}{k}$ and $\tau_2 : \Dsquare{v_2}{w_2}{k}{l}$,
    we have a \textbf{vertical composition}
    \[ \sqvcomp{\tau_1}{\tau_2} : \Dsquare{v_1 \cdot v_2}{w_1 \cdot w_2}{h}{l}; \]
  \item\label{double-cat:horizontal-id-square} for all $v : x \rightarrow y$,
    we have a \textbf{horizontal identity} \[ \sqhid{v} : \Dsquare{v}{v}{\horid{x}}{\horid{y}}; \]
  \item\label{double-cat:horizontal-comp-square} for all squares $\tau_1 : \Dsquare{v_1}{v_2}{h_1}{k_1}$ and $\tau_2 : \Dsquare{v_2}{v_3}{h_2}{k_2}$,
    we have a \textbf{horizontal composition}
    \[ \sqhcomp{\tau_1}{\tau_2} : \Dsquare{v_1}{v_3}{\horcomp{h_1}{h_2}}{\horcomp{k_1}{k_2}}; \]
  \item\label{double-cat:lunitor} for all $h : \hor{x}{y}$, we have a \textbf{left unitor} \[ \Dlunitor{h} : \Dsquare{\id{x}}{\id{y}}{\horcomp{\id{x}}{h}}{h}; \]
  \item\label{double-cat:runitor} for all $h : \hor{x}{y}$, we have a \textbf{right unitor} \[ \Drunitor{h} : \Dsquare{\id{x}}{\id{y}}{\horcomp{h}{\id{y}}}{h}; \]
  \item\label{double-cat:associator} for all $h_1 : \hor{w}{x}$, $h_2 : \hor{x}{y}$, and $h_3 : \hor{y}{z}$, we have an \textbf{associator} \[ \Dassociator{h_1}{h_2}{h_3} : \Dsquare{\id{w}}{\id{z}}{\horcomp{h_1}{(\horcomp{h_2}{h_3})}}{\horcomp{(\horcomp{h_1}{h_2})}{h_3}}. \]
\end{enumerate}
This data is required to satisfy several laws, stating, in particular,
that horizontal identities and horizontal composition are functorial,
and that the left unitor, right unitor, and associator are natural transformations.
In addition, we have the \textbf{triangle} and \textbf{pentagon} law.
Their description can be found in \cref{fig:triangle} and \cref{fig:pentagon}.
\end{definition}

\begin{figure*}
\[
\begin{tikzcd}[column sep=4em, row sep=2.3em]
  x & y & y & z \\
  x & y && z
  \arrow[""{name=0, anchor=center, inner sep=0}, "h", "\shortmid"{marking}, from=1-1, to=1-2]
  \arrow[""{name=1, anchor=center, inner sep=0}, "{\id{y}}", "\shortmid"{marking}, from=1-2, to=1-3]
  \arrow["k", "\shortmid"{marking}, from=1-3, to=1-4]
  \arrow[from=1-1, to=2-1]
  \arrow[""{name=2, anchor=center, inner sep=0}, "h"', "\shortmid"{marking}, from=2-1, to=2-2]
  \arrow[from=1-2, to=2-2]
  \arrow[""{name=3, anchor=center, inner sep=0}, "k"', "\shortmid"{marking}, from=2-2, to=2-4]
  \arrow[from=1-4, to=2-4]
  \arrow["{\Dlunitor{k}}"{description}, shorten <=5pt, shorten >=5pt, Rightarrow, from=1, to=3]
  \arrow["{\sqhid{v}}"{description}, shorten <=4pt, shorten >=4pt, Rightarrow, from=0, to=2]
\end{tikzcd}
\quad
\depEq
\quad
\begin{tikzcd}[column sep=4em, row sep=2.3em]
  x & y & y & z \\
  x & y & y & z \\
  x && y & z
  \arrow[from=1-1, to=2-1]
  \arrow[from=2-1, to=3-1]
  \arrow[""{name=0, anchor=center, inner sep=0}, "h", "\shortmid"{marking}, from=1-1, to=1-2]
  \arrow["{\id{y}}", "\shortmid"{marking}, from=1-2, to=1-3]
  \arrow["k", "\shortmid"{marking}, from=1-3, to=1-4]
  \arrow["{\id{y}}", "\shortmid"{marking}, from=2-2, to=2-3]
  \arrow[""{name=1, anchor=center, inner sep=0}, "k", "\shortmid"{marking}, from=2-3, to=2-4]
  \arrow[from=1-4, to=2-4]
  \arrow[""{name=2, anchor=center, inner sep=0}, "h"', "\shortmid"{marking}, from=3-1, to=3-3]
  \arrow[""{name=3, anchor=center, inner sep=0}, "k"', "\shortmid"{marking}, from=3-3, to=3-4]
  \arrow[from=2-4, to=3-4]
  \arrow[from=2-3, to=3-3]
  \arrow[""{name=4, anchor=center, inner sep=0}, "h", "\shortmid"{marking}, from=2-1, to=2-2]
  \arrow["{\sqhid{k}}"{description}, shorten <=4pt, shorten >=4pt, Rightarrow, from=1, to=3]
  \arrow["{\Drunitor{h}}"{description}, shorten <=5pt, shorten >=5pt, Rightarrow, from=4, to=2]
  \arrow["{\Dassociator{h}{\id{y}}{k}}"{description}, shorten <=13pt, shorten >=13pt, Rightarrow, from=0, to=1]
\end{tikzcd}
\]
\caption{The triangle equation}
\label{fig:triangle}
\end{figure*}

\begin{figure*}
\[
\begin{tikzcd}[column sep=4em, row sep=2.3em]
  v & w & x & z \\
  v & x & y & z \\
  v &&& z
  \arrow[""{name=0, anchor=center, inner sep=0}, "{h_1}", "\shortmid"{marking}, from=1-1, to=1-2]
  \arrow["{h_2}", "\shortmid"{marking}, from=1-2, to=1-3]
  \arrow["{\horcomp{h_3}{h_4}}", "\shortmid"{marking}, from=1-3, to=1-4]
  \arrow[from=1-1, to=2-1]
  \arrow[from=1-4, to=2-4]
  \arrow[""{name=1, anchor=center, inner sep=0}, "{\horcomp{h_1}{h_2}}", "\shortmid"{marking}, from=2-1, to=2-2]
  \arrow["{h_3}", "\shortmid"{marking}, from=2-2, to=2-3]
  \arrow[""{name=2, anchor=center, inner sep=0}, "{h_4}", "\shortmid"{marking}, from=2-3, to=2-4]
  \arrow[from=2-1, to=3-1]
  \arrow[""{name=3, anchor=center, inner sep=0}, "{\horcomp{(\horcomp{(\horcomp{h_1}{h_2})}{h_3})}{h_4}}"', "\shortmid"{marking}, from=3-1, to=3-4]
  \arrow[from=2-4, to=3-4]
  \arrow["{\Dassociator{h_1}{h_2}{\horcomp{h_3}{h_4}}}"{description}, shorten <=13pt, shorten >=13pt, Rightarrow, from=0, to=2]
  \arrow["{\Dassociator{\horcomp{h_1}{h_2}}{h_3}{h_4}}"{description}, shorten <=7pt, shorten >=7pt, Rightarrow, from=1, to=3]
\end{tikzcd}
\quad
\depEq
\quad
\begin{tikzcd}[column sep=4em, row sep=2.3em]
  v & w & x & y & z \\
  v & w & x & y & z \\
  v & w & x & y & z \\
  v & w & x & y & z
  \arrow[""{name=0, anchor=center, inner sep=0}, "{h_1}", "\shortmid"{marking}, from=1-1, to=1-2]
  \arrow[""{name=1, anchor=center, inner sep=0}, "{h_1}"', "\shortmid"{marking}, from=2-1, to=2-2]
  \arrow[from=1-1, to=2-1]
  \arrow[from=1-2, to=2-2]
  \arrow[from=2-1, to=3-1]
  \arrow["v", from=3-1, to=4-1]
  \arrow[""{name=2, anchor=center, inner sep=0}, "{h_1}", "\shortmid"{marking}, from=3-1, to=3-2]
  \arrow[from=1-5, to=2-5]
  \arrow[""{name=3, anchor=center, inner sep=0}, "{h_4}", "\shortmid"{marking}, from=3-4, to=3-5]
  \arrow[from=2-5, to=3-5]
  \arrow[""{name=4, anchor=center, inner sep=0}, "{h_4}"', "\shortmid"{marking}, from=4-4, to=4-5]
  \arrow[from=3-5, to=4-5]
  \arrow[from=3-4, to=4-4]
  \arrow[""{name=5, anchor=center, inner sep=0}, "{h_2}", "\shortmid"{marking}, from=1-2, to=1-3]
  \arrow["{h_3}", "\shortmid"{marking}, from=1-3, to=1-4]
  \arrow["{h_4}", "\shortmid"{marking}, from=1-4, to=1-5]
  \arrow["{h_2}"', "\shortmid"{marking}, from=2-2, to=2-3]
  \arrow["{h_3}"', "\shortmid"{marking}, from=2-3, to=2-4]
  \arrow[""{name=6, anchor=center, inner sep=0}, "{h_4}"', "\shortmid"{marking}, from=2-4, to=2-5]
  \arrow["{h_2}", "\shortmid"{marking}, from=3-2, to=3-3]
  \arrow["{h_3}", "\shortmid"{marking}, from=3-3, to=3-4]
  \arrow["{h_1}"', "\shortmid"{marking}, from=4-1, to=4-2]
  \arrow["{h_2}"', "\shortmid"{marking}, from=4-2, to=4-3]
  \arrow[""{name=7, anchor=center, inner sep=0}, "{h_3}"', "\shortmid"{marking}, from=4-3, to=4-4]
  \arrow["{\sqhid{h_1}}"{description}, shorten <=4pt, shorten >=4pt, Rightarrow, from=0, to=1]
  \arrow["{\sqhid{h_4}}"{description}, shorten <=4pt, shorten >=4pt, Rightarrow, from=3, to=4]
  \arrow["{\Dassociator{h_2}{h_3}{h_4}}"{description}, shorten <=13pt, shorten >=13pt, Rightarrow, from=5, to=6]
  \arrow["{\Dassociator{h_1}{\horcomp{h_2}{h_3}}{h_4}}"{description}, shorten <=19pt, shorten >=19pt, Rightarrow, from=1, to=3]
  \arrow["{\Dassociator{h_1}{h_2}{h_3}}"{description}, shorten <=13pt, shorten >=13pt, Rightarrow, from=2, to=7]
\end{tikzcd}\]
\caption{The pentagon equation}
\label{fig:pentagon}
\end{figure*}

To formulate the laws in \cref{def:double-cat-unfolded}, one needs to use transports.
The necessity of these transports come from the laws of the squares.
For example, if we compose a square $\tau : \Dsquare{v}{w}{h}{k}$ with the identity square,
then we should get the original square $\tau$ back.
However, the square $\sqvcomp{\tau}{\sqvid{k}}$ has different sides than $\tau$,
because the the top and bottom sides of $\sqvcomp{\tau}{\sqvid{k}}$ are composed with identities.
As such, we need the laws for vertical composition in order to state the laws for composition of squares.

\begin{remark}
\label{rem:agda-cats}
Strict double categories can also be defined in an unfolded style.
One can do so by slightly modifying \cref{def:double-cat-unfolded}:
we add the requirement that the horizontal morphisms form a set
and that the unitors and associators are identities.
Such an approach is used in the Agda-categories library \cite{DBLP:conf/cpp/HuC21}.
\end{remark}

However, \cref{def:double-cat-unfolded} is still unsatisfactory for our purposes.
Many notions from double category theory can be derived from the bicategory $\bicatofdoublecats$ of double categories.
For example, equivalences of double categories are the same as adjoint equivalences in $\bicatofdoublecats$ \cite{MR2534210},
monoidal double categories are the same pseudomonoids in $\bicatofdoublecats$ \cite{MR4170469},
and fibrations of double categories are the same as internal Street fibrations \cite{MR4520578}.
For this reason, $\bicatofdoublecats$ plays a prominent role in double category theory.

Since we are working in univalent foundations, we would also like a notion of univalence for double categories
and a univalence principle for them.
This principle can be formulated by saying that $\bicatofdoublecats$ is a univalent bicategory.
All in all, our goals in this paper are
\begin{itemize}
  \item to define the notion of univalent double category;
  \item to define the bicategory $\bicatofdoublecats$ of double categories;
  \item to prove that $\bicatofdoublecats$ is a univalent bicategory.
\end{itemize}

The unfolded definition from \cref{def:double-cat-unfolded} poses several complications for our purposes.
More specifically, proving that $\bicatofdoublecats$ is univalent, would become unfeasible.
This is because we are forced to consider the identity type of double categories, which is rather complicated.
However, by using \emph{displayed bicategories} \cite{DBLP:journals/mscs/AhrensFMVW21}, one can give a simpler proof that $\bicatofdoublecats$ is univalent.
Intuitively, the idea is to break up the definition into smaller layers.
The identity type of each of these layers is simpler, and that simplifies the proof of univalence.

This is the basic philosophy behind the definition of double category that we describe in the remainder of this paper.
More specifically, we take the following steps:
\begin{itemize}
  \item We define the notion of 2-sided displayed categories in \cref{sec:2-sided-displayed}.
    With 2-sided displayed categories, we can describe categories with an additional class of morphisms and squares.
  \item In \cref{sec:bicat-of-double-cat}, we describe the bicategory of double categories.
    We start by defining the displayed bicategory of 2-sided displayed categories,
    and step-by-step we add data and properties to acquire double categories.
    For example, in \cref{def:disp-bicat-horid}, we add horizontal identities to the structure,
    and in \cref{def:disp-bicat-horcomp}, we add a horizontal composition operation.
    Simultaneously, we prove that the resulting bicategory is univalent.
\end{itemize}

Another advantage of our approach
is that we can use it to construct adjoint equivalences and invertible 2-cells of double categories.
We describe this process in \cref{sec:equiv-inv-double}.

\section{2-Sided Displayed Categories}
\label{sec:2-sided-displayed}
The notion of displayed categories was developed by Ahrens and Lumsdaine \cite{DBLP:journals/lmcs/AhrensL19}.
Displayed categories are useful for various purposes, and among those are defining the notion of Grothendieck fibration
and modularly defining univalent categories.
Intuitively, a displayed category represents structure/property of objects and morphisms in some category $\CC$.
Displayed categories consist of a type family of \emph{displayed objects} parametrized by the objects of $\CC$,
and a family of sets of \emph{displayed morphisms} parametrized by the morphisms in $\CC$ and displayed objects.
For example, we have a displayed category of group structures over the category of sets.
The displayed objects over a set $X$ are group structures on $X$,
and the set of displayed morphisms over $f : X \rightarrow Y$ from a group structure $G_X$ over $X$
to a group structure $G_Y$ over $Y$ are proofs that $f$ preserves the group operations.

In this section, we define 2-sided displayed categories --- a variation of the notion of displayed categories.
The difference between 2-sided displayed categories and displayed categories is that displayed categories depend on \emph{one} category,
whereas 2-sided displayed categories depend on \emph{two} categories.
Note that 2-sided displayed categories share many purposes with displayed categories:
they can be used to define univalent categories in a modular way, and they can be used to define 2-sided fibrations \cite{MR4177953,MR0574662}.
However, in this paper we view 2-sided displayed categories in another way, namely as an extra class of morphisms on a category.

\begin{definition}[\coqident{CategoryTheory.TwoSidedDisplayedCats.TwoSidedDispCat}{twosided_disp_cat}]
\label{def:2-sided-disp-cat}
Let $\CC_1$ and $\CC_2$ be categories.
A \textbf{2-sided displayed category} $\CD$ over $\CC_1$ and $\CC_2$ consists of
\begin{enumerate}
  \item for all objects $x_1 : \CC_1$ and $x_2 : \CC_2$ a type $\Dob{\CD}{x_1}{x_2}$ of \textbf{objects over $x$ and $y$}
  \item for all objects $\disp{x} : \Dob{\CD}{x_1}{x_2}$ and $\disp{y} : \Dob{\CD}{y_1}{y_2}$ and morphisms $f_1 : x_1 \rightarrow y_1$ in $\CC_1$ and $f_2 : x_2 \rightarrow y_2$ in $\CC_2$, a set $\Dmor{\disp{x}}{\disp{y}}{f}{g}$ of \textbf{morphisms over $f_1$ and $f_2$}
  \item for every object $\disp{x} : \Dob{\CD}{x_1}{x_2}$ a morphism $\Did{\disp{x}}$ over $\id{x_1}$ and $\id{x_2}$
  \item for all $\disp{f} : \Dmor{\disp{x}}{\disp{y}}{f_1}{f_2}$ and $\disp{g} : \Dmor{\disp{y}}{\disp{z}}{g_1}{g_2}$, a morphism $\disp{f} \cdot \disp{g} : \Dmor{\disp{x}}{\disp{z}}{f_1 \cdot g_1}{f_2 \cdot g_2}$
\end{enumerate} \saveitem
such that the following equations hold.
\begin{enumerate} \restoreitem
  \item \label{2-sided-disp-cat:2-sided-unit} for all $\disp{f} : \Dmor{\disp{x}}{\disp{y}}{f_1}{f_2}$, we have $\disp{f} \cdot \Did{\disp{y}} \depEq \disp{f}$ and $\Did{\disp{x}} \cdot \disp{f} \depEq \disp{f}$;
  \item \label{2-sided-disp-cat:2-sided-assoc} for all $\disp{f} : \Dmor{\disp{w}}{\disp{x}}{f_1}{f_2}$, $\disp{g} : \Dmor{\disp{x}}{\disp{y}}{g_1}{g_2}$, and $\disp{h} : \Dmor{\disp{y}}{\disp{z}}{h_1}{h_2}$, we have $\disp{f} \cdot (\disp{g} \cdot \disp{h}) \depEq (\disp{f} \cdot \disp{g}) \cdot \disp{h}$.
\end{enumerate}
\end{definition}

Here we use the notation $\depEq$ to represent a \emph{dependent} equality, i.e. a path between an element $y_1 : Y(x_1)$ and $y_2 : Y(x_2)$ such that $x_1 = x_2$.
Note that the laws in \cref{2-sided-disp-cat:2-sided-unit,2-sided-disp-cat:2-sided-assoc} in \cref{def:2-sided-disp-cat} are actually dependent equalities.
For examples, if $\disp{f} : \Dmor{\disp{x}}{\disp{y}}{f_1}{f_2}$,
then the left-hand side of $\disp{f} \cdot \Did{\disp{y}} \depEq \disp{f}$ is a morphism that lives over $f_1 \cdot \id{y_1}$ and $f_2 \cdot \id{y_2}$, respectively.
However, the right-hand side lives over $f_1$ and $f_2$,
and thus their types are not equal.
We can solve this by properly using a transport.

Every displayed category $\CD$ over $\CC$ gives rise to a total category $\total{\CD}$ and a functor $\total{\CD} \rightarrow \CC$.
For 2-sided displayed categories, we can do the same.

\begin{definition}[\coqident{CategoryTheory.TwoSidedDisplayedCats.Total}{total_twosided_disp_category}]
\label{def:total}
Let $\CD$ be a 2-sided displayed category over $\CC_1$ and $\CC_2$.
Then we define the \textbf{total category} $\total{\CD}$ to be the category whose objects consists of triples $x_1 : \CC_1$, $x_2 : \CC_2$, and $\disp{x} : \Dob{\CD}{x_1}{x_2}$.
We also define the \textbf{projection functors} $\projl{\CD} : \total{\CD} \rightarrow \CC_1$ and $\projr{\CD} : \total{\CD} \rightarrow \CC_2$
to be the functors that take the first and second coordinate of a triple, respectively.
\end{definition}

Note that every 2-sided displayed category $\CD$ over $\CC_1$ and $\CC_2$ gives rise to a span $\SpanDiag{\CC_1}{\projl{\CD}}{\total{\CD}}{\projr{\CD}}{\CC_2}$ of categories.
Now let us consider some examples of 2-sided displayed categories.

\begin{example}[\coqident{CategoryTheory.TwoSidedDisplayedCats.Examples.Arrow}{arrow_twosided_disp_cat}]
\label{exa:square-disp}
Let $\CC$ be a category.
We define the 2-sided displayed category $\Arr{\CC}$ over $\CC$ and $\CC$ as follows.
\begin{itemize}
  \item The objects over $x$ and $y$ are morphisms $\varphi : x \rightarrow y$.
  \item Suppose that we have morphisms $f : x_1 \rightarrow x_2$, $g : y_1 \rightarrow y_2$, $\varphi_1 : x_1 \rightarrow y_1$, and $\varphi_2 : x_2 \rightarrow y_2$,
    then the set $\Dmor{\varphi_1}{\varphi_2}{f}{g}$ is defined to be the collection of proofs that $f \cdot \varphi_2 = \varphi_1 \cdot g$.
\end{itemize}
The total category $\total{\Arr{\CC}}$ is equivalent to the arrow category of $\CC$.
\end{example}

\begin{example}[\coqident{CategoryTheory.TwoSidedDisplayedCats.Examples.Comma}{comma_twosided_disp_cat}]
\label{exa:comma-disp}
Given functors $F : \CC_1 \rightarrow \CC_3$ and $G : \CC_2 \rightarrow \CC_3$,
we define the 2-sided displayed category $\Comma{F}{G}$ over $\CC_1$ and $\CC_2$:
\begin{itemize}
  \item The objects over $x : \CC_1$ and $y : \CC_2$ are morphisms $\varphi : F(x) \rightarrow G(y)$.
  \item Given morphisms $f : x_1 \rightarrow x_2$, $g : y_1 \rightarrow y_2$, $\varphi_1 : F(x_1) \rightarrow G(y_1)$, and $\varphi_2 : F(x_2) \rightarrow G(y_2)$,
    the set $\Dmor{\varphi_1}{\varphi_2}{f}{g}$ is defined to be the collection of proofs that $F(f) \cdot \varphi_2 = \varphi_1 \cdot G(g)$.
\end{itemize}
The category $\total{\Comma{F}{G}}$ is equivalent to the comma category of $F$ and $G$.
\end{example}

\begin{example}[\coqident{CategoryTheory.TwoSidedDisplayedCats.Examples.Spans}{twosided_disp_cat_of_spans}]
\label{exa:spans-disp}
Let $\CC$ be a category.
We define the 2-sided displayed category $\Span{\CC}$ over $\CC$ and $\CC$:
\begin{itemize}
  \item The objects over $x$ and $y$ are \textbf{spans} from $x : \CC$ to $y : \CC$.
    More concretely, they consist of an object $z : \CC$ and two morphisms $\varphi : z \rightarrow x$ and $\psi : z \rightarrow y$.
  \item Suppose that we have $f : x_1 \rightarrow x_2$ and $g : y_1 \rightarrow y_2$.
    A morphism from $\SpanDiag{x_1}{\varphi_1}{z_1}{\psi_1}{y_1}$ to $\SpanDiag{x_2}{\varphi_2}{z_2}{\psi_2}{y_2}$
    over $f$ and $g$ consists of a morphism $h : z_1 \rightarrow z_2$ such that the following
    diagrams commute.
    \[\begin{tikzcd}
	{x_1} & {z_1} & {y_1} \\
	{x_2} & {z_2} & {y_1}
	\arrow["{\varphi_1}"', from=1-2, to=1-1]
	\arrow["{\varphi_2}", from=2-2, to=2-1]
	\arrow["{\psi_1}", from=1-2, to=1-3]
	\arrow["{\psi_2}"', from=2-2, to=2-3]
	\arrow["f"', from=1-1, to=2-1]
	\arrow["g", from=1-3, to=2-3]
	\arrow["h"{description}, from=1-2, to=2-2]
      \end{tikzcd}\]
\end{itemize}
\end{example}

\begin{example}[\coqident{CategoryTheory.TwoSidedDisplayedCats.Examples.StructuredCospans}{twosided_disp_cat_of_struct_cospans}]
\label{exa:struct-cospan-disp}
Suppose that we have a functor $L : \CC_1 \rightarrow \CC_2$.
We define the 2-sided displayed category $\StructCospan{L}$ over $\CC_1$ and $\CC_1$:
\begin{itemize}
  \item The objects over $x$ and $y$ are \textbf{structured cospans} from $x : \CC_1$ to $y : \CC_1$,
    that is to say, an object $z : \CC_2$ together with morphisms $\Cospan{L(x)}{\varphi}{z}{\psi}{L(y)}$.
  \item Given two structured cospans $\Cospan{L(x_1)}{\varphi_1}{z_1}{\psi_1}{L(y_1)}$ and $\Cospan{L(x_2)}{\varphi_2}{z_2}{\psi_2}{L(y_2)}$,
    and two morphisms $f : x_1 \rightarrow x_2$ and $g : y_1 \rightarrow y_2$,
    a displayed morphism consists of a morphism $h : z_1 \rightarrow z_2$ such that the following diagram commutes
    \[\begin{tikzcd}
	{L(x_1)} & {z_1} & {L(y_1)} \\
	{L(x_2)} & {z_2} & {L(y_1)}
	\arrow["{\varphi_2}"', from=2-1, to=2-2]
	\arrow["{\psi_1}"', from=1-3, to=1-2]
	\arrow["{\psi_2}", from=2-3, to=2-2]
	\arrow["L(f)"', from=1-1, to=2-1]
	\arrow["L(g)", from=1-3, to=2-3]
	\arrow["h"{description}, from=1-2, to=2-2]
	\arrow["{\varphi_1}", from=1-1, to=1-2]
      \end{tikzcd}\]
\end{itemize}
\end{example}

\begin{example}[\coqident{CategoryTheory.TwoSidedDisplayedCats.Examples.Lenses}{twosided_disp_cat_of_lenses}]
\label{exa:lenses-twosided-disp-cat}
Let $\CC$ be a category with chosen binary products.
A \textbf{lens} $l$ from $s$ to $v$ consists of a \textbf{get}-morphism $\lensget{l} : s \rightarrow v$ and a \textbf{put}-morphism $\lensput{l} : v \times s \rightarrow s$ such that
\begin{itemize}
  \item $\lensput{l} \cdot \lensget{l} = \pi_1$;
  \item $\lensget{l} \times \id{s} \cdot \lensput{l} = \id{s}$;
  \item $\id{v} \times \lensput{l} \cdot \lensput{l} = \pi_1 \times (\pi_2 \cdot \pi_2) \cdot \lensput{l}$.
\end{itemize}
Then we define a 2-sided displayed category $\Lens{C}$ over $\CC$ and $\CC$ as follows.
\begin{itemize}
  \item The displayed objects over $s$ and $v$ are lenses from $s$ to $v$.
  \item Given morphisms $f_1 : s_1 \rightarrow s_2$ and $f_2 : v_1 \rightarrow v_2$ and lenses $l_1$ from $s_1$ to $v_1$ and $l_2$ from $s_2$ to $v_2$,
    the displayed morphisms from $l_1$ to $l_2$ over $f_1$ and $f_2$ are proofs that $g_1 \cdot f_2 = f_1 \cdot g_2$ and $p_1 \cdot f_1 = f_2 \times f_1 \cdot p_2$.
\end{itemize}
\end{example}

Our next goal is to define \emph{univalent} 2-sided displayed categories.
To do so, we take the same approach as for categories and for displayed categories.
We first define the notion of isomorphism, and we prove that the identity is an isomorphism.
With that in place, we obtain a map that sends equalities of displayed objects to isomorphisms,
and univalence is formulated by saying that this map is an equivalence of types.

\begin{definition}[\coqident{CategoryTheory.TwoSidedDisplayedCats.Isos}{is_iso_twosided_disp}]
\label{def:disp-iso}
Let $\CD$ be a 2-sided displayed category over $\CC_1$ and $\CC_2$,
and let $f_1 : x_1 \rightarrow y_1$ and $f_2 : x_2 \rightarrow y_2$ be isomorphisms in $\CC_1$ and $\CC_2$ respectively.
In addition, suppose that we have objects $\disp{x} : \Dob{\CD}{x_1}{y_1}$ and $\disp{y} : \Dob{\CD}{x_2}{y_2}$.
Then we say that $\disp{f} : \Dmor{\disp{x}}{\disp{y}}{f_1}{f_2}$ is an \textbf{isomorphism}
if we have a morphism $\disp{f^{-1}} : \Dmor{\disp{y}}{\disp{x}}{\inv{f_1}}{\inv{f_2}}$
such that $\disp{f} \cdot \disp{f^{-1}} \depEq \Did{\disp{x}}$ and $\disp{f^{-1}} \cdot \disp{f} \depEq \Did{\disp{y}}$.
\end{definition}

\begin{proposition}[\coqident{CategoryTheory.TwoSidedDisplayedCats.Isos}{isaprop_is_iso_twosided_disp}]
\label{prop:isaprop-iso}
For every morphism $\disp{f} : \Dmor{\disp{x}}{\disp{y}}{f_1}{f_2}$ over isomorphisms $f_1$ and $f_2$,
the type that $\disp{f}$ is an isomorphism is a proposition.
\end{proposition}

\begin{proposition}[\coqident{CategoryTheory.TwoSidedDisplayedCats.Isos}{id_iso_twosided_disp}]
\label{prop:id-iso}
For all displayed objects $\disp{x} : \Dob{\CD}{x_1}{x_2}$, the identity $\Did{\disp{x}}$ is an isomorphism.
\end{proposition}

\begin{definition}[\coqident{CategoryTheory.TwoSidedDisplayedCats.Univalence}{is_univalent_twosided_disp_cat}]
\label{def:univ}
Let $\CD$ be a 2-sided displayed category over $\CC_1$ and $\CC_2$.
\begin{itemize}
  \item For all objects $x_1 : \CC_1$ and $x_2 : \CC_2$ and displayed objects $\disp{x}, \disp{y} : \Dob{\CD}{x_1}{x_2}$,
    we have a map that sends identities $p : \disp{x} = \disp{y}$ to isomorphisms $\idtoiso{\disp{x}}{\disp{y}}(p) : \Dmor{\disp{x}}{\disp{y}}{\id{x_1}}{\id{x_2}}$.
  \item We say that $\CD$ is \textbf{univalent} if for all $\disp{x}, \disp{y} : \Dob{\CD}{x_1}{x_2}$, the map $\idtoiso{\disp{x}}{\disp{y}}$ is an equivalence of types.
\end{itemize}
\end{definition}

Note that in the formalization, the definition of univalence is equivalent, but formulated slightly differently.
In \cref{def:univ}, we only look at paths $p : \disp{x} = \disp{y}$ between displayed objects lying over the same objects
in the base, whereas in the formalization, we also take paths in the base into account.
Each of the 2-sided displayed categories from \cref{exa:square-disp,exa:comma-disp,exa:spans-disp,exa:struct-cospan-disp,exa:lenses-twosided-disp-cat} is univalent.

\begin{proposition}[{\coqident[is_univalent_total]{CategoryTheory.TwoSidedDisplayedCats.Total}{is_univalent_total_twosided_disp_category}}]
\label{prop:univ-total}
If $\CD$ is a univalent 2-sided displayed category over $\CC_1$ and $\CC_2$, and $\CC_1$ and $\CC_2$ are univalent,
then $\total{\CD}$ is univalent as well.
\end{proposition}

Recall that every 2-sided displayed category gives rise to a span of categories.
Hence, by \cref{prop:univ-total}, every univalent 2-sided displayed category $\CD$ over $\CC_1$ and $\CC_2$ gives rise to a span $\SpanDiag{\CC_1}{\projl{\CD}}{\total{\CD}}{\projr{\CD}}{\CC_2}$ of univalent categories.
To end this section, we define the notions of \emph{2-sided displayed functors} and \emph{2-sided displayed natural transformations}.
These play a prominent role when we define the bicategory of double categories in \cref{sec:bicat-of-double-cat}.

\begin{definition}[\coqident{CategoryTheory.TwoSidedDisplayedCats.DisplayedFunctor}{twosided_disp_functor}]
\label{def:disp-functor}
Suppose that we have 2-sided displayed categories $\CD$ over $\CC_1$ and $\CC_2$, and $\CD'$ over $\CC_3$ and $\CC_4$.
In addition, suppose that we have functors $F_1 : \CC_1 \rightarrow \CC_3$ and $F_2 : \CC_2 \rightarrow \CC_4$.
A \textbf{2-sided displayed functor} $\disp{F}$ over $F_1$ and $F_2$ from $\CD$ to $\CD'$ consists of
\begin{itemize}
  \item a map that assigns to every object $\disp{x} : \Dob{\CD}{x_1}{x_2}$ an object $\disp{F}(\disp{x}) : \Dob{\CD'}{F_1(x_1)}{F_2(x_2)}$;
  \item a map that assigns to every morphism $\disp{f} : \Dmor{\disp{x}}{\disp{y}}{f_1}{f_2}$ a morphism $\disp{F}(\disp{f}) : \Dmor{\disp{F}(\disp{x})}{\disp{F}(\disp{y})}{F_1(f_1)}{F_2(f_2)}$
\end{itemize}
such that $\disp{F}(\Did{\disp{x}}) \depEq \Did{\disp{F}(\disp{x})}$ and $\disp{F} (\disp{f} \cdot \disp{g}) \depEq \disp{F}(\disp{f}) \cdot \disp{F}(\disp{g})$.
\end{definition}

\begin{definition}[\coqident{CategoryTheory.TwoSidedDisplayedCats.DisplayedNatTrans}{twosided_disp_nat_trans}]
\label{def:disp-nat-trans}
Suppose that we have 2-sided displayed categories $\CD$ over $\CC_1$ and $\CC_2$, and $\CD'$ over $\CC_3$ and $\CC_4$.
In addition, suppose that we have functors $F_1, G_1 : \CC_1 \rightarrow \CC_3$ and $F_2, G_2 : \CC_2 \rightarrow \CC_4$,
2-sided displayed functors $\disp{F}$ over $F_1$ and $F_2$ and $\disp{G}$ over $G_1$ and $G_2$,
and natural transformations $\tau_1 : \nattrans{F_1}{G_1}$ and $\tau_2 : \nattrans{F_2}{G_2}$.
A \textbf{2-sided displayed natural transformation} $\disp{\tau}$ over $\tau_1$ and $\tau_2$ from $\disp{F}$ to $\disp{G}$ consists of
a map that assigns to every $\disp{x} : \Dob{\CD}{x_1}{x_2}$ a morphism $\Dmor{\disp{F}(\disp{x})}{\disp{G}(\disp{x})}{\tau_1(x_1)}{\tau_2(x_2)}$
such that the usual naturality condition holds.
\end{definition}

\section{A Recap on (Displayed) Bicategories}
\label{sec:disp-bicat}
Our next goal is to construct the bicategory of double categories.
To do so, we recall in this section the definitions and propositions that we use in the remainder of this paper.
These definitions were originally introduced in \cite{DBLP:journals/mscs/AhrensFMVW21}, and full definitions can be found there.
Recall that a bicategory not only has objects and morphisms, but also 2-cells.
The notion of \emph{displayed bicategory} is similar to that of displayed category.

\begin{definition}[\coqident{Bicategories.DisplayedBicats.DispBicat}{disp_bicat}]
\label{def:disp-bicat}
Let $\CB$ be a bicategory.
A \textbf{displayed bicategory} $\CD$ over $\CB$ consists of
\begin{itemize}
  \item for each object $x : \CB$, a type $\DBob{\CD}{x}$ of objects over $x$;
  \item for all 1-cells $f : x \rightarrow y$ and displayed objects $\disp{x} : \DBob{\CD}{x}$ and $\disp{y} : \DBob{\CD}{y}$,
    a type $\DBmor{\disp{x}}{\disp{y}}{f}$ of 1-cells over $f$;
  \item for all 2-cells $\tau : \cell{f}{g}$ and displayed 1-cells $\disp{f} : \DBmor{\disp{x}}{\disp{y}}{f}$ and $\disp{g} : \DBmor{\disp{x}}{\disp{y}}{g}$, a set $\DBcell{\disp{f}}{\disp{g}}{\tau}$ of 2-cells over $\tau$.
\end{itemize}
In addition, there should be suitable identities, composition, unitors, and associators, and the usual coherence laws should be satisfied.  
\end{definition}

There are numerous examples of displayed bicategories and they are discussed in \cite{DBLP:journals/mscs/AhrensFMVW21},
and we quickly recall the ones that we need in \cref{sec:bicat-of-double-cat}.
If we have displayed bicategories $\CD_1$ and $\CD_2$ over $\CB$,
then we have a displayed bicategory $\proddispbicat{\CD_1}{\CD_2}$ over $\CD$
whose displayed objects, 1-cells, and 2-cells are pairs of displayed objects, 1-cells, and 2-cells of $\CD_1$ and $\CD_2$ respectively.
The full subbicategory can also be defined using a displayed bicategory: if we have a predicate $P$ on the objects of a bicategory $\CB$,
then we define a displayed bicategory over $\CB$ whose displayed objects over $x$ are proofs of $P(x)$,
and whose displayed 1-cells and 2-cells are inhabitants of the unit type.

Every displayed bicategory gives rise to a \emph{total bicategory}.

\begin{definition}[\coqident{Bicategories.DisplayedBicats.DispBicat}{total_bicat}]
\label{def:total-bicat}
Given a displayed bicategory $\CD$ over $\CB$, we define its \textbf{total bicategory} as the bicategory
whose objects are given by pairs of objects $x : \CB$ and $\disp{x} : \DBob{\CD}{x}$.
The 1-cells and 2-cells are defined similarly.
\end{definition}

\emph{Univalent bicategories} are defined in a similar way as univalent categories, but there is a slight difference.
For categories, univalence is expressed by saying that identity of objects is equivalent to
isomorphisms between objects.
For bicategories on the other hand, we formulate univalence in two steps.
First of all, we say that identity of 1-cells is equivalent to invertible 2-cells between them.
This is called \emph{local univalence} in \cite{DBLP:journals/mscs/AhrensFMVW21}.
Secondly, we say that identity of objects is equivalent to adjoint equivalences between them.
In \cite{DBLP:journals/mscs/AhrensFMVW21}, this is called \emph{global univalence}.
Then a \emph{univalent bicategory} is a bicategory that is both locally and globally univalent.
Similarly, we define \emph{univalent displayed bicategories}.
The key theorem for univalent displayed bicategories is the following.

\begin{proposition}[\coqident{Bicategories.DisplayedBicats.DispUnivalence}{total_is_univalent_2}]
Let $\CD$ be a univalent displayed bicategory over a univalent bicategory $\CB$.
Then $\total{\CD}$ is univalent.
\end{proposition}

One key application of univalence for bicategories, is \emph{equivalence induction}.
More specifically, to prove some property for every invertible 2-cell,
it suffices to only consider identity 2-cells.
Similarly, to prove some property for every adjoint equivalence,
one only has to show it for identity equivalences.
This is similar to path induction in homotopy type theory \cite{hottbook,rijke2022introduction}.

\begin{proposition}[\coqident{Bicategories.Core.Univalence}{J_2_0}]
\label{prop:equiv-induction}
Let $\CB$ be a univalent bicategory,
and suppose that for all objects $x, y : \CB$,
we have a predicate $P$ on adjoint equivalences $\adjequiv{x}{y}$.
Then $P$ holds for every adjoint equivalence
if $P$ holds for $\id{x} : \adjequiv{x}{x}$ for every $x : \CB$.
\end{proposition}

\begin{proposition}[\coqident{Bicategories.Core.Univalence}{J_2_1}]
\label{prop:inv2cell-induction}
Let $\CB$ be a univalent bicategory,
and suppose that for all objects $x, y : \CB$ and 1-cells $f, g : x \rightarrow y$,
we have a predicate $P$ on invertible 2-cells $\invcell{f}{g}$.
Then $P$ holds for every invertible 2-cell
if $P$ holds for $\id{f} : \invcell{f}{f}$ for every $f : x \rightarrow y$.
\end{proposition}

\section{The Bicategory of Double Categories}
\label{sec:bicat-of-double-cat}
In this section, we define the bicategory of univalent double categories,
and we prove that this bicategory is univalent.
The notion of \emph{displayed bicategory} plays a key role in this construction \cite{DBLP:journals/mscs/AhrensFMVW21}.

The construction proceeds in several steps.
We start in \cref{def:disp-bicat-of-twosided-disp-cat} by defining a displayed bicategory $\dispbicatoftwosideddispcat$ over the bicategory $\bicatofcats$ of univalent categories,
and the objects over $\CC$ are 2-sided displayed categories $\CD$ over $\CC$ and $\CC$.
If we look at the total bicategory $\bicatoftwosideddispcat$ of this displayed bicategory,
then the objects consists of a category $\CC$ and a 2-sided displayed $\CD$ over $\CC$ and $\CC$.
This means that we have a category with an extra class of morphisms and a class of squares.

To obtain a the bicategory of double categories,
we need to add more structure.
We define two displayed bicategories $\dispbicathorid$ and $\dispbicathorcomp$ over $\bicatoftwosideddispcat$ in \cref{def:disp-bicat-horid} and \cref{def:disp-bicat-horcomp}.
The displayed bicategory $\dispbicathorid$ adds horizontal identities to the structure, and $\dispbicathorcomp$ adds horizontal compositions.
By taking their product and the total bicategory, we obtain the bicategory $\bicathoridcomp$, of which the objects consists of a category, horizontal morphisms, squares, horizontal identities, and compositions.

Next we define displayed bicategories $\dispbicatlunitor$, $\dispbicatrunitor$, and $\dispbicatassociator$ over $\bicathoridcomp$.
These add the left unitor, the right unitor, and the associator to the structure.
Again we take their product and the total bicategory to obtain the bicategory $\bicatunitorassociator$.
Finally, we define $\bicatofdoublecats$ as a full subbicategory of $\bicatunitorassociator$:
the predicate we use, expresses the triangle and pentagon coherence.

At each step, we prove that the relevant displayed bicategories are univalent.
The machinery of displayed bicategories allows us to combine all of this to conclude that $\bicatofdoublecats$ is univalent.
The advantage of using displayed bicategories over a direct approach is that the proof of univalence becomes simpler and more modular.
This is because the displayed approach allows us to consider the identity of each part individually,
and we are able to reuse results (e.g., the bicategory of univalent categories is univalent).

The main idea behind this construction is that we can split up the definition of a double category into several layers.
Instead of looking at the whole, we look at these layers separately, and that allows for reusability and modularity.
This is also why the notion of 2-sided displayed category plays an important role in this construction:
it is one of the layers to define double categories.

\begin{definition}[\coqident{Bicategories.DisplayedBicats.Examples.DispBicatOfTwoSidedDispCat}{disp_bicat_twosided_disp_cat}]
\label{def:disp-bicat-of-twosided-disp-cat}
The displayed bicategory $\dispbicatoftwosideddispcat$ over $\bicatofcats$ is defined as follows:
\begin{itemize}
  \item The displayed objects over $\CC$ are univalent 2-sided displayed categories $\CD$ over $\CC$ and $\CC$.
  \item The displayed morphisms from $\CD_1$ to $\CD_2$ over $F : \CC_1 \rightarrow \CC_2$ are 2-sided displayed functors $\disp{F}$ over $F$ and $F$ from $\CD_1$ to $\CD_2$.
  \item The displayed 2-cells from $\disp{F}$ to $\disp{G}$ over $\tau : \nattrans{F}{G}$ are 2-sided displayed natural transformations over $\tau$ and $\tau$ from $\disp{F}$ to $\disp{G}$.
\end{itemize}
We define $\bicatoftwosideddispcat$ to be $\total{\dispbicatoftwosideddispcat}$.
\end{definition}

An object of $\bicatoftwosideddispcat$ consists of a univalent category $\CC$
and a univalent 2-sided displayed category $\CD$ over $\CC$ and $\CC$.
If we compare this to \cref{def:double-cat-unfolded},
then we already got the data from \cref{double-cat:vertical-cat,double-cat:hor-mor,double-cat:squares,double-cat:vertical-id-square,double-cat:vertical-comp-square}.
The vertical category is given by $\CC$,
the horizontal morphisms from $x$ to $y$ are given by the displayed objects $\Dob{\CD}{x}{y}$,
and the squares $\Dsquare{v_1}{v_2}{h}{k}$ are given by displayed morphisms $\Dmor{h}{k}{v_1}{v_2}$.
The vertical identity and composition for squares is given by the identity and composition in $\CD$, respectively,
and similarly for the laws involving vertical composition of squares.

\begin{proposition}[{\coqident[univalent_2_twosided_disp_cat]{Bicategories.DisplayedBicats.Examples.DispBicatOfTwoSidedDispCat}{disp_univalent_2_disp_bicat_twosided_disp_cat}}]
\label{prop:disp-bicat-of-twosided-disp-cat-univ}
The displayed bicategory $\dispbicatoftwosideddispcat$ is univalent.
\end{proposition}

\subsection{Identities and Composition}
\label{sec:id-and-comp}
Next we add horizontal identities (\cref{double-cat:hor-id,double-cat:horizontal-id-square} in \cref{def:double-cat-unfolded}
and composition (\cref{double-cat:hor-comp,double-cat:horizontal-comp-square} in \cref{def:double-cat-unfolded}),
in the form of two displayed bicategories over $\bicatoftwosideddispcat$.
To define the first one, we define when a 2-sided displayed category supports \emph{horizontal identities}.

\begin{definition}[\coqident{Bicategories.DoubleCategories.DoubleCategoryBasics}{hor_id}]
\label{def:disp-hor-id}
Let $\CC$ be a category and let $\CD$ be a 2-sided displayed category over $\CC$ and $\CC$.
Then we say that $\CD$ has \textbf{horizontal identities} if
\begin{enumerate}
  \item for all $x : \CC$, we have a displayed object $\doublecathorid{\CD}{x} : \Dob{\CD}{x}{x}$;
  \item for all morphisms $v : x \rightarrow y$, we have a displayed morphism $\doublecathoridsq{\CD}{v} : \Dmor{\doublecathorid{\CD}{x}}{\doublecathorid{\CD}{y}}{v}{v}$;
\end{enumerate}
such that $\doublecathoridsq{\CD}{\id{x}} = \id{\doublecathorid{\CD}{x}}$
and $\doublecathoridsq{\CD}{v_1 \cdot v_2} = \doublecathoridsq{\CD}{v_1} \cdot \doublecathoridsq{\CD}{v_2}$.
\end{definition}

We also define when a 2-sided displayed functor \emph{preserves horizontal identities}.

\begin{definition}[\coqident{Bicategories.DoubleCategories.DoubleFunctors.Basics}{double_functor_hor_id}]
\label{def:disp-pres-hor-id}
Let $\CD$ be a 2-sided displayed category over $\CC$ and $\CC$ and let $\CD'$ be a 2-sided displayed category over $\CC'$ and $\CC'$.
Suppose that we have a functor $F : \CC \rightarrow \CC'$ and a 2-sided displayed functor $\disp{F}$ from $\CD$ to $\CD'$ over $F$ and $F$,
and that $\CD$ and $\CD'$ have horizontal identities.
Then we say that $\disp{F}$ \textbf{preserves horizontal identities} if for all $x : \CC_1$
we have a natural square $\doublefunctorhorid{F}{x} : \Dmor{\doublecathorid{\CD'}{F(x)}}{F(\doublecathorid{\CD}{x})}{\id{F(x)}}{\id{F(x)}}$.
\end{definition}

The precise naturality condition for the square can be found in the formalization.
In addition, note that we consider \emph{lax} double functors: we do not require $\doublefunctorhorid{F}{x}$ to be invertible.

\begin{definition}[\coqident{Bicategories.DoubleCategories.BicatOfDoubleCats}{disp_bicat_twosided_disp_cat_hor_id}]
\label{def:disp-bicat-horid}
We define the displayed bicategory $\dispbicathorid$ over $\bicatoftwosideddispcat$ as follows:
\begin{itemize}
  \item the displayed objects over a pair of a univalent category $\CC$ and a univalent 2-sided displayed category $\CD$ are horizontal identities for $\CD$ (\cref{def:disp-hor-id});
  \item the displayed 1-cells over a functor $F : \CC_1 \rightarrow \CC_2$ and 2-sided displayed functor $\disp{F}$ from $\CD_1$ to $\CD_2$ that preserve horizontal identities (\cref{def:disp-pres-hor-id});
  \item the displayed 2-cells over a natural transformations $\tau : \nattrans{F}{G}$ and a 2-sided displayed natural transformation $\disp{\tau}$ are proofs that $\disp{\tau}$ preserves horizontal identities. The precise formulation can be found in the formalization.
\end{itemize}
\end{definition}

Next we look at horizontal compositions.

\begin{definition}[\coqident{Bicategories.DoubleCategories.DoubleCategoryBasics}{hor_comp}]
\label{def:disp-hor-comp}
Let $\CC$ be a category and let $\CD$ be a 2-sided displayed category over $\CC$ and $\CC$.
Then we say that $\CD$ has \textbf{horizontal composition} if
\begin{itemize}
  \item for all $h : \Dob{\CD}{x}{y}$ and $k : \Dob{\CD}{y}{z}$, we have a displayed object $\doublecathorcomp{h}{k} : \Dob{\CD}{x}{z}$;
  \item for all displayed morphisms $s_1 : \Dmor{h_1}{h_2}{v_1}{v_2}$ and $s_2 : \Dmor{k_1}{k_2}{v_2}{v_3}$, we have a displayed morphism $\sqhcomp{s_1}{s_2} : \Dmor{\doublecathorcomp{h_1}{k_1}}{\doublecathorcomp{h_2}{k_2}}{v_1}{v_3}$;
\end{itemize}
such that
\begin{itemize}
  \item $\sqhcomp{\id{h}}{\id{k}} = \id{\doublecathorcomp{h}{k}}$.
  \item $\sqhcomp{(s_1 \cdot t_1)}{(s_2 \cdot t_2)} = (\sqhcomp{s_1}{s_2}) \cdot (\sqhcomp{t_1}{t_2})$.
\end{itemize}
\end{definition}

\begin{definition}[\coqident{Bicategories.DoubleCategories.DoubleFunctors.Basics}{double_functor_hor_comp}]
\label{def:disp-pres-hor-comp}
Let $\CD$ be a 2-sided displayed categories over $\CC$ and $\CC$ and let $\CD'$ be 2-sided displayed categories over $\CC'$ and $\CC'$.
Suppose that we have a functor $F : \CC \rightarrow \CC'$ and a 2-sided displayed functor $\disp{F}$ from $\CD$ to $\CD'$ over $F$ and $F$,
and that $\CD_1$ and $\CD_2$ have horizontal identities.
Then we say that $\disp{F}$ \textbf{preserves horizontal compositions} if for all $h : \Dob{\CD}{x}{y}$ and $k : \Dob{\CD}{y}{k}$
we have a natural square $\doublefunctorhorcomp{F}{h}{k} : \Dsquare{\doublecathorcomp{F(h)}{F(k)}}{F(\doublecathorcomp{h}{k})}{\id{F(x)}}{\id{F(z)}}$.
\end{definition}

\begin{definition}[\coqident{Bicategories.DoubleCategories.BicatOfDoubleCats}{disp_bicat_twosided_disp_cat_hor_comp}]
\label{def:disp-bicat-horcomp}
The displayed bicategory $\dispbicathorcomp$ over $\bicatoftwosideddispcat$ is defined as follows:
\begin{itemize}
  \item the displayed objects over a pair of a univalent category $\CC$ and a univalent 2-sided displayed category $\CD$ are horizontal composition for $\CD$ (\cref{def:disp-hor-comp});
  \item the displayed 1-cells over a functor $F : \CC_1 \rightarrow \CC_2$ and 2-sided displayed functor $\disp{F}$ from $\CD_1$ to $\CD_2$ that preserve horizontal composition (\cref{def:disp-pres-hor-comp});
  \item the displayed 2-cells over a natural transformations $\tau : \nattrans{F}{G}$ and a 2-sided displayed natural transformation $\disp{\tau}$ are proofs that $\disp{\tau}$ preserves horizontal composition. The precise formulation can be found in the formalization.
\end{itemize}
\end{definition}

We define $\dispbicathoridcomp$ to be $\proddispbicat{\dispbicathorid}{\dispbicathorcomp}$,
and we define $\bicathoridcomp$ to be the total bicategory of $\dispbicathoridcomp$.

\begin{proposition}[\coqident{Bicategories.DoubleCategories.BicatOfDoubleCats}{disp_univalent_2_disp_bicat_twosided_disp_cat_id_hor_comp}]
\label{prop:disp-bicat-horidcomp-univ}
The displayed bicategory $\dispbicathoridcomp$ is univalent.
\end{proposition}

\subsection{Unitors and Associators}
At this point, we obtained the bicategory $\bicathoridcomp$,
and the objects of that bicategory consists of a univalent category $\CC$,
a univalent displayed $\CD$ over $\CC$ and $\CC$,
together with horizontal identities (\cref{def:disp-hor-id})
and horizontal compositions (\cref{def:disp-hor-comp}).
This corresponds to \cref{double-cat:vertical-cat,double-cat:hor-mor,double-cat:hor-id,double-cat:hor-comp,double-cat:squares,double-cat:vertical-id-square,double-cat:vertical-comp-square,double-cat:horizontal-id-square,double-cat:horizontal-comp-square} in \cref{def:double-cat-unfolded}
and now we look at \cref{double-cat:lunitor,double-cat:runitor,double-cat:associator} from \cref{def:double-cat-unfolded}.
For each of these, we define a displayed bicategory over $\bicathoridcomp$.
Due to space constraints, we only say how the displayed objects of those displayed bicategories are defined.

\begin{definition}[\coqident{Bicategories.DoubleCategories.BicatOfDoubleCats}{disp_bicat_lunitor}]
\label{def:disp-bicat-lunitor}
We define a displayed bicategory $\dispbicatlunitor$ over $\bicathoridcomp$
whose displayed objects over a univalent category $\CC$ and a univalent 2-sided displayed category $\CD$ with horizontal identities and compositions
consists of a natural isomorphism $\Dmor{\doublecathorcomp{\doublecathorid{\CD}{x}}{h}}{h}{\id{x}}{\id{y}}$ for each $h : \Dob{\CD}{x}{y}$.
\end{definition}

\begin{definition}[\coqident{Bicategories.DoubleCategories.BicatOfDoubleCats}{disp_bicat_runitor}]
\label{def:disp-bicat-runitor}
We define a displayed bicategory $\dispbicatrunitor$ over $\bicathoridcomp$
whose displayed objects over a univalent category $\CC$ and a univalent 2-sided displayed category $\CD$ with horizontal identities and compositions
consists of a natural isomorphism $\Dmor{\doublecathorcomp{h}{\doublecathorid{\CD}{y}}}{h}{\id{x}}{\id{y}}$ for each $h : \Dob{\CD}{x}{y}$.
\end{definition}

\begin{definition}[\coqident{Bicategories.DoubleCategories.BicatOfDoubleCats}{disp_bicat_lassociator}]
\label{def:disp-bicat-associator}
We define a displayed bicategory $\dispbicatassociator$ over $\bicathoridcomp$
whose displayed objects over a univalent category $\CC$ and a univalent 2-sided displayed category $\CD$ with horizontal identities and compositions
consists of a natural isomorphism
\[
\Dmor{\doublecathorcomp{h_1}{(\doublecathorcomp{h_2}{h_3})}}{\doublecathorcomp{(\doublecathorcomp{h_1}{h_2})}{h_3}}{\id{w}}{\id{z}}
\]
for all $h_1 : \Dob{\CD}{w}{x}$, $h_2 : \Dob{\CD}{x}{y}$, and $h_3 : \Dob{\CD}{y}{z}$.
\end{definition}

We define $\dispbicatunitorassociator$ to be $\proddispbicat{\dispbicatlunitor}{\proddispbicat{\dispbicatrunitor}{\dispbicatassociator}}$,
and we define $\bicatunitorassociator$ to be the total bicategory of $\dispbicatunitorassociator$.

\begin{proposition}[\coqident{Bicategories.DoubleCategories.BicatOfDoubleCats}{is_univalent_2_bicat_unitors_and_associator}]
\label{prop:disp-bicat-unitor-associator-univ}
The displayed bicategory $\dispbicatunitorassociator$ is univalent.
\end{proposition}

\begin{definition}[\coqident{Bicategories.DoubleCategories.BicatOfDoubleCats}{bicat_of_double_cats}]
\label{def:bicat-of-double-cat}
We define the bicategory $\bicatofdoublecats$ of double categories as the full subbicategory $\fullsubbicat{\bicatunitorassociator}{}$
where the predicate expresses that the triangle and pentagon laws are satisfied.
\end{definition}

\begin{theorem}[\coqident{Bicategories.DoubleCategories.BicatOfDoubleCats}{is_univalent_2_bicat_of_double_cats}]
\label{thm:univalent-bicat-of-double-cat}
The bicategory $\bicatofdoublecats$ is univalent.
\end{theorem}

The proof of \cref{thm:univalent-bicat-of-double-cat} uses Voevodsky's univalence axiom.

The objects of $\bicatofdoublecats$ collect all the data and properties mentioned in this section.
Each of these data and properties correspond to some part of \cref{def:double-cat-unfolded}.
For that reason, we can conclude that these two definitions of double categories are actually equivalent.

\begin{theorem}[\coqident{Bicategories.DoubleCategories.DoubleCatsEquivalentDefinitions}{double_cat_weq_univalent_doublecategory}]
\label{thm:double-cat-def-equiv}
The type of objects of $\bicatofdoublecats$ is equivalent to the type of double categories as defined in \cref{def:double-cat-unfolded}.
\end{theorem}

\begin{remark}
	Notice the unfolded definition in \cref{def:double-cat-unfolded} does not include a univalence condition, as we characterize univalent double categories as the ones for which the corresponding $2$-sided displayed category under the equivalence given in \cref{thm:double-cat-def-equiv} is univalent in the sense of \cref{def:univ}.
\end{remark}

The 1-cells in $\bicatofdoublecats$ are \emph{lax double functors}.
They consist of an underlying functor and 2-sided displayed functor that preserve horizontal identities and compositions
as described in \cref{def:disp-pres-hor-id,def:disp-pres-hor-comp}.
Finally, every 2-cell in $\bicatofdoublecats$ has an underlying natural transformation and 2-sided displayed natural transformation.

\begin{remark}
\label{rem:univ-double-cat}
Note that the double categories in $\bicatofdoublecats$ are univalent,
and this univalence condition means that the underlying category and 2-sided
displayed category are univalent.
From this, we see that objects in $\bicatofdoublecats$ are the same as pseudocategories
internal to the bicategory of univalent categories.

In \cite[Example~9.3]{up}, a notion of univalent double bicategory is defined such that identities correspond to \emph{gregarious equivalences} of double bicategories. In the particular case where the underlying bicategory given by objects, vertical morphisms and squares is a category (meaning the 1-morphisms form a set and the assignment from identities of 1-morphisms to 2-morphisms is an equivalence), our notion of univalence coincides with the notion introduced in \cite{up}.
\end{remark}

\section{Examples of Double Categories}
\label{sec:exa-double-cats}
Now we construct several examples of double categories using \cref{def:bicat-of-double-cat}.
All of the double categories considered here are univalent.

\begin{example}[\coqident{Bicategories.DoubleCategories.Examples.SquareDoubleCat}{square_double_cat}]
\label{exa:square-double-cat}
Let $\CC$ be a univalent category.
In \cref{exa:square-disp}, we defined the 2-sided displayed category $\Arr{\CC}$ over $\CC$ and $\CC$.
This gives rise to a double category as follows.
\begin{itemize}
  \item the horizontal identities are given by the identity morphism;
  \item horizontal composition is given by the composition of morphisms.
\end{itemize}
The unitality and associativity of horizontal composition reduce to the ordinary laws of composition for morphisms.
All laws involving squares hold because the type of morphisms in a category is a set.
\end{example}

\begin{example}[\coqident{Bicategories.DoubleCategories.Examples.KleisliDoubleCat}{kleisli_double_cat}]
\label{exa:kleisli-double-cat}
Let $T$ be a monad on a univalent category $\CC$.
In \cref{exa:comma-disp}, we defined the 2-sided displayed category $\Comma{F}{G}$ for arbitrary functors $F : \CC_1 \rightarrow \CC_3$
and $G : \CC_2 \rightarrow \CC_3$.
We take $F$ to be the identity on $\CC$ and $G$ to be the endofunctor underlying $T$.
Concretely, we look at $\Comma{\id{\CC}}{T}$, meaning that the horizontal morphisms are morphisms in the Kleisli category of $T$.
We obtain the following double category.
\begin{itemize}
  \item the horizontal identities are given by the unit of $T$;
  \item given morphisms $h : x \rightarrow T(y)$ and $k : y \rightarrow T(z)$, their horizontal composition is defined as the following composition
    \[\begin{tikzcd}
	x & {T(y)} & {T(T(z))} & {T(z)}
	\arrow["h", from=1-1, to=1-2]
	\arrow["{T k}", from=1-2, to=1-3]
	\arrow["{\mu_z}", from=1-3, to=1-4]
      \end{tikzcd}\]
\end{itemize}
The construction of the unitors and associators for this double category reduces to proving unitality and associativity of composition in the Kleisli category.
\end{example}

One way to instantiate \cref{exa:kleisli-double-cat}
is by taking $\CC$ to be $\SET$
and $T$ to be the power set monad.
Note that morphisms from $X$ to $Y$ in the Kleisli category of the power set monad
are the same as relations between $X$ and $Y$.
Hence, the resulting double category has functions as vertical morphisms,
and relations as horizontal morphisms.

Note that both \cref{exa:square-double-cat,exa:kleisli-double-cat}
are strict double categories.
In both cases, the type of horizontal morphisms is a set and unitality
and associativity for horizontal composition holds up to equality.

\begin{example}[\coqident{Bicategories.DoubleCategories.Examples.SpansDoubleCat}{spans_double_cat}]
\label{exa:span-double-cat}
Let $\CC$ be a univalent category with pullbacks.
We defined the 2-sided displayed category $\Span{\CC}$ in \cref{exa:spans-disp}.
This gives rise to a double category.
\begin{itemize}
  \item The horizontal identity on an object $x : \CC$ is given by the span $\SpanDiag{x}{\id{x}}{x}{\id{x}}{x}$.
  \item Suppose that we have spans $\SpanDiag{x_1}{\varphi_1}{z_1}{\psi_1}{y_1}$ and $\SpanDiag{x_2}{\varphi_2}{z_2}{\psi_2}{y_2}$.
    Their composition is given by the following span
    
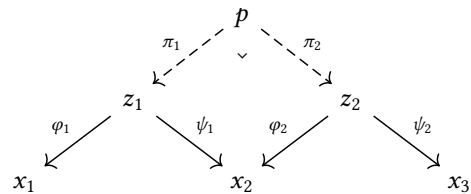
\[\begin{tikzcd}
	&& p \\
	& {z_1} && {z_2} \\
	{x_1} && {x_2} && {x_3}
	\arrow["{\varphi_1}"', from=2-2, to=3-1]
	\arrow["{\psi_1}", from=2-2, to=3-3]
	\arrow["{\varphi_2}"', from=2-4, to=3-3]
	\arrow["{\psi_2}", from=2-4, to=3-5]
	\arrow["{\pi_1}"', dashed, from=1-3, to=2-2]
	\arrow["{\pi_2}", dashed, from=1-3, to=2-4]
	\arrow["\lrcorner"{anchor=center, pos=0.125, rotate=-45}, draw=none, from=1-3, to=3-3]
      \end{tikzcd}\]
    Here $p$ is the pullback of $\psi_1$ and $\varphi_2$.
\end{itemize}

To construct the right unitor of this double category, we consider the following diagram
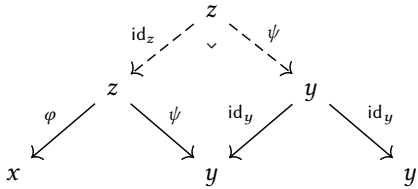
\[\begin{tikzcd}
	&& z \\
	& z && {y} \\
	{x} && {y} && {y}
	\arrow["\varphi"', from=2-2, to=3-1]
	\arrow["\psi", from=2-2, to=3-3]
	\arrow["{\id{y}}"', from=2-4, to=3-3]
	\arrow["{\id{y}}", from=2-4, to=3-5]
	\arrow["{\id{z}}"', dashed, from=1-3, to=2-2]
	\arrow["\psi", dashed, from=1-3, to=2-4]
	\arrow["\lrcorner"{anchor=center, pos=0.125, rotate=-45}, draw=none, from=1-3, to=3-3]
\end{tikzcd}\]
The square in this diagram is a pullback, and from this, we get the desired isomorphism for the right unitor.
Similarly, we can define the left unitor and the associator.
The proofs of the triangle and pentagon laws follow by diagram chasing,
and details can be found in the formalization.
\end{example}

\begin{example}[\coqident{Bicategories.DoubleCategories.Examples.StructuredCospansDoubleCat}{structured_cospans_double_cat}]
\label{exa:struct-cospan-double-cat}
Suppose that we have a functor $L : \CC_1 \rightarrow \CC_2$ between univalent categories and suppose that $\CC_2$ has pushouts.
In \cref{exa:struct-cospan-disp} we defined the 2-sided displayed category $\StructCospan{L}$ over $\CC_1$ and $\CC_1$ of structured cospans.
This gives rise to a double category.
\begin{itemize}
  \item The horizontal identity on an object $x : \CC_1$ is given by the cospan \[\Cospan{L(x)}{\id{L(x)}}{L(x)}{\id{L(x)}}{L(x)}.\]
  \item The construction of the horizontal composition is dual to how horizontal composition is defined in \cref{exa:span-double-cat}.
\end{itemize}
\end{example}

\begin{example}[\cite{riley2018categories}, \coqident{Bicategories.DoubleCategories.Examples.LensesDoubleCat}{lenses_double_cat}]
Let $\CC$ be a univalent category with chosen binary products.
Then we define the double category of lenses of $\CC$ as follows.
In \cref{exa:lenses-twosided-disp-cat}, we defined the 2-sided displayed category $\Lens{\CC}$.
We construct a double category from it as follows.
\begin{itemize}
  \item The identity lens from $x : \CC$ to $x : \CC$ is given by $\id{x} : x \rightarrow x$ and $\pi_1 : x \times x \rightarrow x$.
  \item Suppose we have lenses $l_1$ from $x$ to $y$ and $l_2$ from $y$ to $z$.
    Then we have a lens $\horcomp{l_1}{l_2}$ from $x$ to $z$ such that $\lensget{\horcomp{l_1}{l_2}} = \lensget{l_1} \cdot \lensget{l_2}$,
    and such that $\lensput{\horcomp{l_1}{l_2}}$ is the following composition.
    \[\begin{tikzcd}[column sep=2.3em]
	{z \times x}
        &
        {(z \times x)\times x}
        &&&
        {y \times x}
        &
        x
	\arrow["{\lensput{l_1}}", from=1-5, to=1-6]
	\arrow["{\langle \idI , \pi_2 \rangle}", from=1-1, to=1-2]
	\arrow["{((\idI \times \lensget{l_1}) \cdot \lensput{l_2}) \times \idI}", from=1-2, to=1-5]
      \end{tikzcd}\]
\end{itemize}
\end{example}

Note that there are different gadgets called ``lenses'' in the literature.
The lenses by Clarke~\cite[Def.~3.20]{Clarke2023} are, more specifically, ``delta lenses''.
The double category of delta lenses has, as objects, (small) categories, as horizontal morphisms functors between categories, and vertical morphisms delta-lenses, that is, functors equipped with an extra ``lifting operation'' --- see \cite[Def.~2.1]{Clarke2023} for details.
Squares are suitable commutative squares of functors.

\section{Equivalences of Double Categories}
\label{sec:equiv-inv-double}
In this section, we give sufficient conditions to show that a 1-cell in $\bicatofdoublecats$ is an adjoint equivalence (\cref{thm:adj-equiv-double-cat}),
and that a 2-cell in $\bicatofdoublecats$ is invertible (\cref{thm:inv2cell-double-cat}).
Since these proofs are similar, we only discuss how \cref{thm:adj-equiv-double-cat} is proven.
Let us first give conditions for when a 2-cell in $\bicatofdoublecats$ is invertible.

\begin{theorem}[\coqident{Bicategories.DoubleCategories.InvertiblesAndEquivalences}{invertible_double_nat_trans_weq}]
\label{thm:inv2cell-double-cat}
Let $\tau$ be a 2-cell in $\bicatofdoublecats$.
Then $\tau$ is an invertible 2-cell if and only if its underlying natural transformation and 2-sided displayed natural transformation are pointwise isomorphisms.
\end{theorem}

To characterize adjoint equivalences, we need the notion of a \emph{strong double functor}.

\begin{definition}[\coqident{Bicategories.DoubleCategories.DoubleCats}{is_strong_double_functor}]
\label{def:strong-double-functor}
Let $F$ be a lax double functor.
We say that $F$ is a \textbf{strong double functor} if
$\doublefunctorhorid{F}{x}$ and $\doublefunctorhorcomp{F}{h}{k}$ are isomorphisms for all suitably typed $x$, $h$, and $k$.
\end{definition}

\begin{theorem}[{\coqident[adjoint_equivalence_double_functor_weq]{Bicategories.DoubleCategories.InvertiblesAndEquivalences}{left_adjoint_equivalence_lax_double_functor_weq}}]
\label{thm:adj-equiv-double-cat}
Let $L : \CC_1 \rightarrow \CC_2$ be a 1-cell $\bicatofdoublecats$.
Then $L$ is an adjoint equivalence if and only if $L$ is a strong double functor
and $L$ is an adjoint equivalence in $\bicatoftwosideddispcat$.
\end{theorem}

We give a sketch of our proof of \cref{thm:adj-equiv-double-cat};
it follows the construction of $\bicatoftwosideddispcat$ in \cref{sec:bicat-of-double-cat}.
As such, we first show that $L$ lifts to an adjoint equivalence in $\bicathoridcomp$,
then we show that $L$ lifts to an adjoint equivalence in $\bicatunitorassociator$,
and finally, we conclude that $L$ gives rise to an adjoint equivalence in $\bicatofdoublecats$.

Next we show that $L$ lifts to an adjoint equivalence in $\bicathoridcomp$,
and to do so, we construct a displayed adjoint equivalence over $L$ in both $\dispbicathorid$ and $\dispbicathorcomp$.
We simplify this construction by using induction over adjoint equivalences (\cref{prop:equiv-induction})
for which we use that the bicategory $\dispbicathorcomp$ is univalent.
Intuitively, this allows us to assume that $L$ is the identity equivalence.
More concretely, we show the following.

\begin{lemma}
\label{lem:adj-equiv-double-cat-id}
Let $\CD : \bicatoftwosideddispcat$.
Suppose that $I_1$ and $I_2$ are objects over $\CD$ in $\dispbicathorid$,
and that $f : \DBmor{I_1}{I_2}{\id{\CD}}$.
Note that $f$ consists of a natural square $\tau : \Dsquare{\id{x}}{\id{x}}{\doublecathorid{\CD'}{x}}{\doublecathorid{\CD}{x}}$ for each $x$.
Then $f$ is a displayed adjoint equivalence
if
$\tau(x)$ is an isomorphism for every $x$.
\end{lemma}

In our situation, the assumption in \cref{lem:adj-equiv-double-cat-id} follows from the fact that $L$ preserves the identity up to isomorphism.
Similarly, we can construct a displayed adjoint equivalence over $L$ in $\dispbicathorcomp$,
and this gives us the adjoint equivalence in $\bicathoridcomp$.

To lift $L$ to an adjoint equivalence in $\bicatunitorassociator$,
we need to construct displayed adjoint equivalences over $L$ in $\dispbicatlunitor$, $\dispbicatrunitor$, and $\dispbicatassociator$.
Note that each of these displayed bicategories live over $\bicathoridcomp$.
Again we use \cref{prop:equiv-induction}, so we assume that $L$ is the identity.
Constructing the displayed adjoint equivalences then follows from diagram chasing,
and the precise proof can be found in the formalization.

To conclude \cref{thm:adj-equiv-double-cat},
we note that $\bicatofdoublecats$ is defined as a full subbicategory of $\bicatunitorassociator$.
Since adjoint equivalences in full subbicategories of some bicategory $\CB$ are the same as adjoint equivalences in $\CB$,
we get the desired adjoint equivalence in $\bicatofdoublecats$.

For the converse, we first note that whenever $L$ is an adjoint equivalence $\bicatofdoublecats$,
then $L$ is an adjoint equivalence in $\bicatoftwosideddispcat$.
This is because pseudofunctors preserve adjoint equivalence.
To show that $L$ is a strong double functor, we use \cref{prop:equiv-induction},
so it suffices to show that the identity is a strong double functor.
This follows from the fact that the identity is an isomorphism.

Note that Shulman proves \cref{thm:adj-equiv-double-cat} for framed bicategories in a different way \cite[Corollary 7.9]{MR2534210}.
Whereas our proof follows the construction of $\bicatofdoublecats$ via displayed bicategories
and makes use of induction over adjoint equivalences,
Shulman's proof makes use of fully faithful and essentially surjective strong double functors.

\section{Formalizing Weak Double Categories}
\label{sec:upcoming}
In \cref{sec:bicat-of-double-cat} we constructed the univalent bicategory of univalent double categories, establishing that identities in the type of objects capture the data of equivalences of double categories. This is a direct continuation of \cite{DBLP:journals/mscs/AhrensFMVW21}, where the authors constructed the univalent bicategory of univalent categories. However, while the notion of equivalence of categories is established and consistent throughout the literature, there are several relevant notions of equivalences of double categories, resulting in several relevant univalence principles. 

In this work we focused on equivalences coming from pseudo-double category theory, known as \emph{vertical equivalences}, following conventions in \cite[Section 3.5.5]{grandis2020doublecat} or \citep[Proposition 12.3.21]{johnsonyau2021doublecat}. In this approach the definition of equivalence depends on a choice of direction in which the double categorical data is assumed to be strict (in our case the vertical direction, whereas, for example in \cite[Section 3.5.5]{grandis2020doublecat}, the horizontal direction). 

On the other side, there is another important notion of equivalence of double categories in the literature, known as \emph{gregarious equivalences of double categories} \cite[Slide 16]{campbell2020gregarious}, which is symmetric and does not come from categorical equivalences in either direction\footnote{This differs from gregarious equivalences in a double category as discussed in \cref{rem:univ-double-cat}. The clashing terminology stems from the fact that gregarious equivalences in a double category are an important ingredient towards defining gregarious equivalences between double categories.}. This necessitates a formalization of double categorical structures that are non-strict in both the horizontal and vertical direction.

In future work, we define and then formalize weak double categories (\cref{def:weakdoublecat}), building on Verity's work on \emph{double bicategories} \cite[Definition 1.4.1]{MR2844536}. We then build on our understanding of the univalence principle \cite{up} to define and study univalent weak double categories (\cref{def:univalent weak double cat}) and establish a univalence principle for gregarious equivalences \cite[Slide 20]{rasekh2023itpdouble}. 
This illustrates how the structural features and the chosen strictness of the formalization directly correlates with its equivalences and univalence principle. We can already witness this phenomenon in \cite{DBLP:journals/mscs/AhrensFMVW21}, where the authors formalize $2$-categories and bicategories with diverging univalence principles.

Concretely we will focus on the following definitions and conjectures. 

\begin{definition} \label{def:weakdoublecat}
	A \emph{pre-pseudo-category} is a pseudo-category object in the $(3,1)$-category of $2$-groupoids. A \emph{pseudo-category} is a pre-pseudo category whose hom-types are groupoids. A \emph{weak double category} is a pseudo-category object in pseudo-categories with a set of squares. 
\end{definition}

Weak double categories defined in this way are not identical to double bicategories, However, we conjecture that they are related.

\begin{conjecture}
	There is a fully faithful embedding of weak double categories into Verity's double bicategories with essential image given by double bicategories in which the given horizontal (vertical) bicategory coincides with the bicategory given by objects, horizontal (vertical) morphisms and squares whose vertical (horizontal) edges are identities.
\end{conjecture}

The condition that the objects in the essential image satisfy should be thought of as a \emph{bicategorical completeness condition}, which guarantees that we have well-defined underlying horizontal and vertical bicategories. This condition has already been considered in \cite[Example 9.3]{up}, where a double bicategory satisfying such a condition is called a \emph{doubly weak double category}. Assuming this conjecture, we can summarize the relation between (strict) double categories, weak double categories and double bicategories via the following diagram of inclusions:
\[
\begin{tikzcd}[column sep=1.2cm]
	\parbox{1.1cm}{\noindent Double  \vspace{-0.1cm} \\ Category }   
	\arrow[r, hookrightarrow] & 
	\parbox{2cm}{\noindent Weak Double \vspace{-0.1cm} \\ Category } 
	\arrow[r, hookrightarrow] & 
	\parbox{1.1cm}{\noindent Double \vspace{-0.1cm} \\ Bicategory } 
\end{tikzcd}
\]
As discussed above, we can build on this first conjecture to define \emph{univalent weak double categories}. 

\begin{definition} \label{def:univalent weak double cat}
	A \emph{univalent weak double category} is a weak double category such that its corresponding double bicategory is univalent.
\end{definition}
 
Unwinding definitions, an identity in the type of objects corresponds to a gregarious equivalence in a double category (as discussed in \cref{rem:univ-double-cat}).

\begin{remark} \label{rem:symmetric}
	A key aspect of this definition of a weak double category and univalence that distinguishes it from the notion of double category studied throughout this paper, is that its definition is symmetric and not dependent on choosing a direction.
\end{remark} 

The difference between univalent double categories (in the sense of \cref{def:bicat-of-double-cat}) and their weak version introduced here is indeed non-trivial, as can be seen by the following conjectural example, which should be understood as the categorified version of \cref{exa:span-double-cat}.

\begin{conjecture}
	(Enriched) categories, (enriched) functors and (enriched) profunctors assemble into a univalent weak double category, but not into a univalent double category.
\end{conjecture}

Weak double categories assemble into a tricategory, with a sub-tricategory consisting of univalent weak double categories. We now have the following major conjecture, analogous to \cref{thm:univalent-bicat-of-double-cat}.  

\begin{conjecture}
	The tricategory of univalent weak double categories is univalent.
\end{conjecture}

\section{Conclusion}
\label{sec:conclusion}
In this paper, we constructed the univalent bicategory of univalent double categories.
The main tool in the construction is the notion of 2-sided displayed categories,
which represent categories with an extra class of morphisms and squares.
We also characterized the adjoint equivalences and invertible 2-cells in the bicategory of univalent double categories,
and in that characterization, we made use of univalence at several points.
Finally, we gave numerous examples of univalent double categories.
Among our examples are the double categories of lenses and of structured cospans.

There are numerous ways to extend on this work.
An interesting special case of double categories is given by \emph{framed bicategories} \cite{MR2534210}.
We can obtain a univalent bicategory of univalent framed bicategories by extending the work in \cref{sec:bicat-of-double-cat}:
we take a full subbicategory of $\bicatofdoublecats$ that expresses that the the double category is framed (i.e., some functor is a fibration).
However, currently framed bicategories are not considered in our formalization.
Furthermore, in many applications, one would like to have more structure on a double category,
such as a (symmetric) monoidal structure.
Such structures can conveniently be defined by looking at pseudomonoids in $\bicatofdoublecats$.
To construct a univalent bicategory of (symmetric) monoidal double categories,
one would need to combine ideas from \cite{DBLP:journals/mscs/AhrensFMVW21,DBLP:journals/lmcs/VeltriW21} and \cite{DBLP:conf/types/WullaertMA22}.

In addition, our notion of univalent double category is unable to capture univalent categories with profunctors.
This is because we do not have a category of univalent categories, but only a bicategory.
This is another situation where the right solution is to pursue a formalization of double bicategories \cite{MR2844536} and its suitable notion of univalence  \cite[Example 9.3]{up}.
\begin{acks}
  We gratefully acknowledge the work by the Coq development team in providing the Coq proof assistant and surrounding infrastructure, as well as their support in keeping UniMath compatible with Coq.
  We are very grateful to Mike Shulman for answering our questions about profunctors.
\end{acks}

\bibliographystyle{ACM-Reference-Format}
\bibliography{literature.bib}

\end{document}